\documentstyle[11pt,newlfont,amscd,amssymb]{amsart}
\newcommand{\nc}{\newcommand}

\nc{\one}{\mbox{\bf 1}}
\nc{\invtensor}{\underset{\leftarrow}{\otimes}}
\nc{\rlarrows}{\begin{picture}(1,0.4)
                \put(0,-0.1){\makebox(1,0.2){$\leftarrow$}}
                \put(0,0.1){\makebox(1,0.2){$\ra$}}
              \end{picture}}
\nc{\rra}{\begin{picture}(1,0.4)
                \put(0,-0.1){\makebox(1,0.2){$\lra$}}
                \put(0,0.1){\makebox(1,0.2){$\lra$}}
              \end{picture}}
\nc{\Left}{\mathbf L}  
\nc{\Right}{\mathbf R} 
\nc{\gr}{\operatorname{gr}}
\nc{\Ho}{\operatorname{Ho}}
\nc{\alt}{\operatorname{alt}}
\nc{\Sym}{\operatorname{Sym}}
\nc{\sym}{\operatorname{sym}}
\nc{\id}{\operatorname{id}}
\nc{\Der}{\operatorname{Der}}
\nc{\im}{\operatorname{Im}}
\nc{\Ker}{\operatorname{Ker}}
\nc{\coker}{\operatorname{Coker}}
\nc{\Col}{\operatorname{Col}}
\nc{\ter}{\operatorname{ter}}
\nc{\intl}{\operatorname{int}}
\nc{\out}{\operatorname{out}}

\nc{\TN}{{\cal N}}
\nc{\Norm}{\operatorname{N}}
\nc{\Nor}{\operatorname{N}}
\nc{\Tor}{\operatorname{Tor}}
\nc{\res}{\operatorname{res}}
\nc{\Stab}{\operatorname{Stab}}
\nc{\Hom}{\operatorname{Hom}}
\nc{\chom}{\CH\!o\!m}
\nc{\uhom}{\CH\!o\!m}
\nc{\End}{\operatorname{End}}
\nc{\holim}{\operatorname{holim}}
\nc{\dirlim}{\underset{\rightarrow}{\lim}\,}
\nc{\invlim}{\underset{\leftarrow}{\lim}\,}
\nc{\CB}{\operatorname{\bf CB}}
\nc{\com}{\operatorname{co}}
\nc{\Tot}{\operatorname{Tot}}
\nc{\Th}{\operatorname{Th}}
\nc{\Cech}{\check{C}}
\nc{\Spec}{\operatorname{Spec}}
\nc{\Spf}{\operatorname{Spf}}
\nc{\MC}{\operatorname{MC}}
\nc{\U}{\operatorname{U}}
\nc{\Diff}{{\cal D}\mbox{\em iff}}
\nc{\Mor} {{\cal M}or}
\nc{\Ob}{\operatorname{Ob}}
\nc{\cone}{\widehat}
\nc{\Coder}{\operatorname{Coder}}
\nc{\pr}{\operatorname{pr}}
\nc{\diag}{\operatorname{diag}}

\nc{\Mod}{{\mathtt{mod}}}       
\nc{\Modf}{{\mathtt{modf}}}       
\nc{\Modg}{{\mathtt{modg}}}       
\nc{\Ab}{{\mathtt {Ab}}}          
\nc{\Alg}{{\mathtt {Alg}}} 
\nc{\Algf}{{\mathtt {Algf}}} 
\nc{\Algg}{{\mathtt {Algg}}} 
\nc{\Coalg}{{\mathtt {Coalg}}} 
         
\nc{\dgc}{{\mathtt{dgc}}}
\nc{\dgca}{{\mathtt{dgca}}}
\nc{\dgcu}{{\mathtt{dgcu}}}
\nc{\dgcuf}{{\mathtt{dgcuf}}}
\nc{\dgcf}{{\mathtt{dgcf}}}
\nc{\dgcg}{{\mathtt{dgcg}}}
\nc{\dgcc}{{\mathtt{dgccc}}}

\nc{\dgl}{{\mathtt{dglie}}}
\nc{\dgla}{{\mathtt{dgla}}}
\nc{\dglf}{{\mathtt{dglf}}}
\nc{\dglg}{{\mathtt{dglg}}}
\nc{\dga}{{\mathtt{dga}}}
\nc{\art}{{\mathtt {art}}}
\nc{\dgar}{{\mathtt {dgart}^{\leq 0}}}
\nc{\simpl}{\Delta^0\Ens}
\nc{\Coll}{{\mathtt{Coll}}}
\nc{\Kan}{{\mathtt {Kan}}}

\nc{\Grp}{{\mathtt {Grp}}}
\nc{\Cat}{{\mathtt {Cat}}}
\nc{\Ens}{{\mathtt {Ens}}}
\nc{\op}{{\operatorname{op}}}
\nc{\Op}{{\mathtt{Op}}}

\nc{\Lie}{{\mathtt{LIE}}}
\nc{\Com}{{\mathtt{COM}}}

\nc{\pa}{\partial}

\nc{\CA}{\cal A}
\nc{\CF}{\cal F}
\nc{\CG}{\cal G}
\nc{\CH}{\cal H}
\nc{\CI}{\cal I}
\nc{\CJ}{\cal J}
\nc{\CO}{\cal O}
\nc{\CP}{\cal P}
\nc{\CU}{\cal U}
\nc{\CC}{\cal C}
\nc{\CDD}{\cal D}
\nc{\CL}{\cal L}
\nc{\CM}{\cal M}
\nc{\CS}{\cal S}
\nc{\CT}{\cal T}

\nc{\fa}{\frak a}
\nc{\fg}{\frak g}
\nc{\fk}{\frak k}
\nc{\fh}{\frak h}
\nc{\fm}{\frak m}
\nc{\fn}{\frak n}
\nc{\fS}{\frak S}
\nc{\fI}{\frak I}
\nc{\fA}{\frak A}

\nc{\nen}{\newenvironment}
\nc{\ol}{\overline}
\nc{\ul}{\underline}
\nc{\ra}{\rightarrow}
\nc{\lra}{\longrightarrow}
\nc{\lla}{\longleftarrow}
\nc{\Lra}{\Longrightarrow}
\nc{\Lla}{\Longleftarrow}
\nc{\Llra}{\Longleftrightarrow}
\nc{\hra}{\hookrightarrow}
\nc{\iso}{\overset{\sim}{\lra}}

\nc{\Thm}[1]{Theorem~\ref{#1}}
\nc{\Prop}[1]{Proposition~\ref{#1}}
\nc{\Lem}[1]{Lemma~\ref{#1}}
\nc{\Cor}[1]{Corollary~\ref{#1}}
\nc{\Conj}[1]{Conjecture~\ref{#1}}
\nc{\Claim}[1]{Claim~\ref{#1}}
\nc{\Defn}[1]{Definition~\ref{#1}}
\nc{\Exa}[1]{Example~\ref{#1}}
\nc{\Rem}[1]{Remark~\ref{#1}}
\nc{\Note}[1]{Note~\ref{#1}}

\nen{thm}[1]{\label{#1}{\bf Theorem.\ } \em}{}
\nen{prop}[1]{\label{#1}{\bf Proposition.\ } \em}{}
\nen{lem}[1]{\label{#1}{\bf Lemma.\ } \em}{}
\nen{klem}[1]{\label{#1}{\bf Key Lemma.\ } \em}{}
\nen{cor}[1]{\label{#1}{\bf Corollary.\ } \em}{}
\nen{conj}[1]{\label{#1}{\bf Conjecture.\ } \em}{}
\nen{claim}[1]{\label{#1}{\bf Claim.\ } \em}{}

\nen{defn}[1]{\label{#1}{\bf Definition.\ } }{}
\nen{exa}[1]{\label{#1}{\bf Example.\ } }{}

\nen{rem}[1]{\label{#1}{\em Remark.\ } }{}
\nen{note}[1]{\label{#1}{\em Note.\ } }{}
\nen{exer}[1]{\label{#1}{\em Exercise.\ } }{}

\begin{document}

\title[]{DG coalgebras as formal stacks}
\author{Vladimir Hinich}
\address{Dept. of Mathematics and Computer Science, University of Haifa,
Mount Carmel, Haifa 31905 Israel}

\maketitle

\section{Introduction}

\subsection{}
In this paper we provide the category of unital (dg, unbounded) coalgebras
$\dgcu(k)$ over a field $k$ of characteristic zero with a structure 
of  a simplicial closed model category (SCMC=simplicial CMC) --- 
see ~\ref{scmc}.

This structure generalizes the one defined by Quillen~\cite{rht} in 1969
for 2-reduced unital coalgebras. The major difference is that our
notion of weak equivalence is strictly stronger (see the example 
in~\ref{counter}) that that of quasi-isomorphism.

The structure of the proof is close to that of Quillen --- we use the pair
of adjoint functors connecting $\dgcu(k)$ with the category  $\dgl(k)$
of Lie dg algebras over $k$ and the existing SCMC structure on $\dgl(k)$
--- see~\cite{haha}. The proof is technically more difficult since
the spectral sequence arguments do not work here and one should carefully 
work in the tensor category of filtered complexes instead.

The main theorems~\ref{cmc-uc} and ~\ref{2} are proven in Sections 3 -- 7.

In the rest of the paper (and of the Introduction) we try to persuade
the reader that unital dg coalgebras provide a proper context to talk 
about formal deformations in characteristic zero. 

In this approach formal deformation problems over a field $k$ of 
characteristic zero are described (up to an equivalence) by functors
\begin{equation}
 F:\dgar(k)\to\simpl
\label{functor}
\end{equation}
from the category of non-positively graded artinian local dg $k$-algebras
to the category of simplicial sets. These functors are represented
(in a homotopy category sense) by unital coalgebras from $\dgcu(k)$
--- see~\ref{nerve.1}.
Such a  coalgebra plays role of the coalgebra of distributions concentrated
at a point and it is defined uniquely up to a weak equivalence in $\dgcu(k)$.

In Section 8 we study properties of the functors~(\ref{functor}) 
which appear as {\em nerves} of dg Lie algebras. In Section 9 we calculate
some elementary examples. In particular, we identify in~\ref{moduli-G} 
the coalgebra of the formal stack of deformations of a principal $G$-bundle $P$
on a scheme $X$ with the standard complex of the Lie algebra 
$\Right\Gamma(X,\fg_P)$. Note that this 
answer is well-known (see, for instance,~\cite{ka} for a similar problem);
however it was not so clear (at least for the author)  what did this answer 
describe.  

\subsection{Tangent Lie algebra}
We suggest to interpret a unital dg coalgebra as the coalgebra of 
distributions concentrated at a point of a would-be-a-space.

\subsubsection{}
The Quillen functor $\CL:\dgcu(k)\to\dgl(k)$ to the category of dg Lie 
algebras (see~\cite{rht}, App.~B or~\ref{quillen} below) is then interpreted
as the tangent Lie algebra functor. \Thm{2} claims that the adjoint pair
$(\CL,\CC)$ of functors establishes an equivalence between the homotopy
categories of $\dgcu(k)$ and of $\dgl(k)$. This means that a formal dg stack
can be (uniquely up to weak equivalence) reconstructed from the
homotopy type of its tangent Lie algebra.

\subsubsection{} 
Formal schemes define, of course, coalgebras concentrated in 
degree zero. The corresponding tangent Lie algebra is concentrated in 
strictly positive degrees. More generally, a formal stack $X\in\dgcu(k)$
satisfying $H^i(\CL(X))=0$ for $i\leq 0$, is called {\em a formal space}.
Equivalently, this means that $X$ is weakly equivalent to a coalgebra
concentrated in non-negative degrees. Formal spaces have also a description
in terms of the functor on Artin rings they represent --- 
see~\ref{formal.spaces.represent} below. 
 
\subsection{Functors on Artin rings}
Thus, we consider unital coalgebras as ``the most general'' formal
(dg) stacks concentrated in a  point. It is reasonable to 
describe them as functors on formal spaces --- as one defines stacks
using functors on affine schemes. One can go even further and take into 
account that any formal space is a filtered colimit of finite dimensional ones
which now take form $A^*$ where $A\in\dgar(k)$.

\subsubsection{}
Any formal stack $X\in\dgcu(k)$ gives rise to a {\em deformation
functor}
$$ \widetilde{X}:\dgar(k)\to\simpl$$
which is defined up to homotopy equivalence. This corresponds to the usual 
description of stacks as 2-functors from affine schemes to groupoids. 
The deformation functor enjoys nice exactness properties
(see~\ref{properties}). Given the tangent Lie algebra $\fg=\CL(X)$,
the functor $ \widetilde{X}$ can be described as  the {\em nerve} 
$\Sigma_{\fg}$ of the dg Lie algebra $\fg$ (see~\cite{ddg} and~\ref{def.nerve})
defined by the formula
$$\Sigma_{\fg}((A,\fm))_n=\MC(\fm\otimes\Omega_n\otimes\fg)$$
for $(A,\fm)\in\dgar(k)$, where $\MC(\_)$ denotes the collection of 
Maurer-Cartan elements of a dg Lie algebra and $\Omega_n$ is the algebra 
of polynomial differential forms on the standard $n$-simplex.

The nerve of a dg Lie algebra is homotopy equivalent to the
Deligne groupoid (cf.~\cite{gm1},~\cite{ddg}) if 
$(\fm\otimes\fg)^i=0\text{ for }i<0$. 

\subsubsection{}
\label{formal.spaces.represent}
If $X$ is a formal space, the restriction of $\widetilde{X}$ to 
the category $\art(k)$ of artinian $k$-algebras concentrated in degree zero,
is a functor to discrete simplicial sets (i.e., essentially, to $\Ens$)
--- see~\ref{f.space}. 

This means in particular, that the restriction of $\widetilde{X}$ to $\art(k)$
is representable in usual sense by $H^0(X)$.
 
\subsection{Rational spaces}

A very ``non-geometric'' class of unital coalgebras is the Quillen's 
category $\dgcu_2(\Bbb{Q})$ of 2-reduced unital coalgebras which is one of the
models for simply connected rational homotopy types.

One can easily calculate the deformation functor defined by a
simply connected rational homotopy type. Let $\fg\in\dgl(\Bbb{Q})$
be the Lie algebra model for it. This is the tangent Lie algebra
of the corresponding unital coalgebra. One has $\fg^i=0$ for $i\geq 0$.

Let $(A,\fm)\in\dgar(\Bbb{Q})$. 
Then $\Sigma_{\fg}(A)$ is a simply-connected rational space
and its homotopy type corresponds to the Lie algebra
$\fm\otimes\fg$ --- see~\ref{rht-example}.

\subsection{Coarse moduli}
Let $\fg\in\dgl(k)$. One might be willing to consider the functor 
$\Pi_{\fg}:(A,\fm)\mapsto \pi_0(\Sigma_{\fg}(A))$
on the category $\art(k)$ as the 
``coarse moduli'' space for the deformation problem defined by $\fg$. 
Usually this functor is not representable (except for the case described in
\ref{formal.spaces.represent}). 
However, it admits a hull in the sense of~\cite{sc}
\footnote{translated to the language of coalgebras} 
which can be easily constructed using a dg Lie subalgebra $\fh$ which is a 
$1$-truncation of 
$\fg$ as in~\cite{gm2}. A $1$-truncation $\fh$ being chosen, the hull
of the functor $\Pi_{\fg}$ can be described as $H^0(\CC(\fh))$ --- 
see~\ref{coarse}.

The choice of the truncation $\fh$ is not unique though the resulting
coalgebra is unique up to a non-canonical isomorphism by a general result 
of~\cite{sc}. One might ask whether the 1-truncation $\fh$ of $\fg$ is 
unique up to a quasi-isomorphism. 
This is obviously so in the case of~\cite{gm2} where $H^i(\fg)=0$ for 
$i\leq 0$. We doubt this is true in general.

\subsection{}Two general ideas

{\bf

Idea 1: Any reasonable formal deformation problem in characteristic zero 
can be described by Maurer-Cartan elements of an appropriate dg Lie algebra

Idea 2: Moduli spaces should admit a natural sheaf of dg commutative algebras 
as a structure sheaf
}

have been spelled out by different people during the last years 
(Drinfeld~\cite{d}, Feigin, Deligne, Kontsevich~\cite{ko}). In this paper
we tried to show that these two claims are essentially equivalent to the one
saying that any reasonable formal deformation problem can be described
by a representable functor on dg artinian rings with values in simplicial 
sets.

\subsection{Acknowledgements} A part of this work was written during
my stay at IHES. I express my gratitude to the Institute for hospitality.

I am very grateful to O.~Gabber who explained to me why formal smoothness
is a local property in fpqc topology.

\section{Preliminaries}
In this Section we fix some standard notations and definitions.
 
Throughout this Section $k$ is a commutative ring containing $\Bbb Q$.

\subsection{Unital coalgebras}

Let $X$ be a cocommutative dg coalgebra  $X$ over $k$ with comultiplication
$\Delta:X\to X\otimes X$ and counit $\epsilon:X\to k$. 

Recall that an element $1\in X$ is called {\em a group-like element} if

(1) $d(1)=0;$ (2) $\Delta(1)=1\otimes 1;$ (3) $\epsilon(1)=1$.

A choice of a group-like element $1\in X$ defines a decomposition
$$ X=k\cdot 1\oplus\ol{X}$$
where $\ol{X}=\ker(\epsilon)$. This defines an increasing  filtration on $X$
by the formula
\begin{equation}
X_n=\ker(X\overset{\Delta_n}{\lra}X^{\otimes n+1}\to\ol{X}^{\otimes n+1})
\label{filtration}
\end{equation}
where $\Delta_n$ is the $n$-th iteration of $\Delta$.
\subsubsection{}
\begin{defn}{uc}(see~\cite{hdtc}) A pair $(X,1)$ consisting of a dg 
cocommutative coalgebra $X$ and a group-like element $1$ is called 
{\em a unital coalgebra} if the filtration~(\ref{filtration}) is exhausting.
\end{defn}

The group-like element defining a unital algebra $X$ is called {\em the unit
of } $X$. The filtration~(\ref{filtration}) of a unital coalgebra is called
{\em the canonical filtration}.

Note that if $k$ is a field then unital coalgebras are just connected 
coalgebras of Quillen (see~\cite{rht}, App. B). In this case the unit
element is unique and it is preserved by any coalgebra map.

\subsubsection{}
The category of unital coalgebras over $k$ will be denoted by $\dgcu(k)$.
The morphisms in it are supposed to preserve the units (this is automatically
fulfilled when $k$ is a field).

\subsubsection{}
Let $\Omega$ be a commutative dg algebra over $k$. Similarly to the above, 
one defines $\Omega$-coalgebras as cocommutative coalgebras in the category
of $\Omega$-modules. Furthermore,
one defines unital $\Omega$-coalgebras as pairs $(X,1)$ with 
a $\Omega$-coalgebra $X$ and
a group-like element $1$ of $X$ such that the filtration~(\ref{filtration})
is exhausting. The category of unital $\Omega$-coalgebras is denoted 
$\dgcu(\Omega)$.

\subsection{Quillen functors} 
\label{quillen}
Recall the definition of the couple of adjoint
functors
\begin{equation}
   \CL:\dgcu(k)\rlarrows\dgl(k):\CC
\label{LandC}
\end{equation}
defined by Quillen in~\cite{rht}, App. B.

\subsubsection{}
Let $X\in\dgcu(k)$, $\ol{X}=\ker{\epsilon}$ in the standard notation.
The dg Lie algebra $\CL(X)$ is defined as follows. As a graded Lie algebras, 
this is the free Lie algebra $F(\ol{X}[-1])$. The differential in $\CL(X)$
is the sum of two parts: the one generated by the differential of $\ol{X}[-1]$,
and the second defined to be the only derivation of the free Lie algebra
$F(\ol{X}[-1])$ whose restriction on $\ol{X}[-1])$ is given by the map
$$\Delta-1\otimes\id-\id\otimes 1:\ol{X}\to\ol{X}\otimes\ol{X}.$$

\subsubsection{}
Let $\fg\in\dgl(k)$. The unital coalgebra $\CC(\fg)$ is defined as follows.
As a graded  coalgebra, this is the cofree cocommutative coalgebra 
$S(\fg[1])$. The differential in $\CC(\fg)$ is the sum of two parts: the one
generated by the differential in $\fg[1]$, and the second defined by its
$1$-component given by the Lie bracket $[,]:\wedge^2\fg\to\fg$.

\subsubsection{}Let $\fg\in\dgl(k)$. Recall that an element $x\in\fg^1$
is called a Maurer-Cartan element if $dx+\frac{1}{2}[x,x]=0$. The set
of Maurer-Cartan elements of $\fg$ is denoted by $\MC(\fg)$.

\subsubsection{} If $X\in\dgcu(k), \fg\in\dgl(k)$, then the complex
$\Hom(\ol{X},\fg)$ admits a natural structure of a dg Lie algebra given
by the formula
$$ [f,g]=\mu (f\otimes g)\Delta$$
where $\Delta$ is the comultiplication in $X$ and $\mu$ is the bracket in 
$\fg$. In particular the set $\MC(X,\fg):=\MC(\Hom(\ol{X},\fg))$ is defined. 
This is the set of {\em twisting functions} from $X$ to $\fg$. 

\subsubsection{}
\begin{thm}{LadjC}The functors $\CL$ and $\CC$ ~(\ref{LandC}) are adoint.
More precisely, for $X\in\dgcu(k), \fg\in\dgl(k)$ one has natural bijections
$$ \Hom(\CL(X),\fg)=\Hom(X,\CC(\fg))=\MC(X,\fg).$$
\end{thm}

See ~\cite{rht}, App. B6.

\subsection{Simplicial closed model categories}
\label{scmc}

We use the axioms (CM1)--(CM5) of~\cite{rht} for the definition of closed
model category (CMC). 

A simplicial category $\CC$ is a collection of objects $\Ob\CC$ together
with a collection of simplicial sets $\uhom(X,Y)$ assigned to each pair
$(X,Y)$ of objects, with strictly associative compositions.

A simplicial category $\CC$ with a CMC structure will
be called a simplicial CMC (or SCMC) if the following Quillen's axiom 
(SM7) ~\cite{ha} is fulfilled

(SM7) Let $i:A\to B$ be a cofibration and $p:X\to Y$ be a fibration in $\CC$.
Then the map of simplicial sets
\begin{equation}
 \uhom(B,X)\to\uhom(A,X)\times_{\uhom(A,Y)} \uhom(B,X)
\label{sm7}
\end{equation}
is a Kan fibration. If, moreover, either $i$ or $p$ is a weak equivalence,
then~(\ref{sm7}) is an acyclic Kan fibration.

{\em Note that we do not include Quillen's axiom (SM0)
claiming the existence of cylinder and path objects  --- see~\cite{ha} ---
in our definition of simplicial closed model category.
}

\subsection{Models for dg Lie algebras}

Recall~\cite{haha} that the category $\dgl(k)$
of dg Lie algebras over a commutative ring $k\supseteq \Bbb Q$ admits a
simplicial model structure. More precisely, one has the following

\begin{thm}{cmc-lie}(see~\cite{haha}, 4.1.1 and 4.8 )
The category $\dgl(k)$ admits a simplicial CMC  structure 
with surjective maps as fibrations and quasi-isomorphisms as weak equivalences.
The simplicial structure on $\dgl(k)$ is defined by the formula
\begin{equation}
\uhom_n(\fg,\fh)=\Hom_{\dgl(k)}(\fg,\Omega_n\otimes\fh).
\label{simp.dgl}
\end{equation}
where $\Omega_n$ is the algebra
of polynomial differential forms on the standard $n$-simplex.
\end{thm}

\subsection{Operad notations}
\label{operads}
Throughout the paper we will sometimes use the language of operads ---
see e.g. ~\cite{hla}. In particular, $\Com$ denotes the operad governing
commutative algebras, $\Lie$ the one for Lie algebras and $\CS$ denotes
the standard Lie operad of~\cite{hla} governing ``strongly homotopy Lie
algebras''. If $\CO$ is an operad and $A$ is
an $\CO$-algebra then $U(\CO,A)$ denotes the corresponding enveloping
algebra. 

We will use sometimes different base tensor categories. If $\CC$ is a
tensor (= symmetric monoidal) category, $\Op(\CC)$ denotes the category
of operads over $\CC$. 

\section{Main theorems}
\label{mainthms}

Now we are ready to formulate the main results of the paper.

\subsection{}
\begin{thm}{cmc-uc}
The category $\dgcu(k)$ of unital coalgebras over a field $k$ of 
characteristic zero admits a simplicial CMC structure. Cofibrations in it
are just injective maps and weak equivalences are the maps $f$ in $\dgcu(k))$
such  that $\CL(f)$ is a quasi-isomorphism.  The simplicial structure
on $\dgcu(k)$ is given by the condition
\begin{equation}
\uhom_n(X,Y)=
\Hom_{\dgcu(\Omega_n)}(\Omega_n\otimes X, \Omega_n\otimes Y)
\label{simp.dgcu}
\end{equation}
where as in~(\ref{simp.dgl}) $\Omega_n$ is the algebra of polynomial 
differential forms on the standard $n$-simplex.
\end{thm}

Note that the property of a morphism $f$ of $\dgcu(k)$ to be a weak 
equivalence is strictly stronger then that of being a quasi-isomorphism --- 
see the counter-example in~\ref{counter}.

\subsection{}
\begin{thm}{2}The adjoint functors $\CL,\CC$ induce an equivalence of the
corresponding homotopy categories
$$\Left\CL:\Ho(\dgcu(k))\rlarrows\Ho(\dgl(k)):\Right\CC.$$
\end{thm}

The proof of the theorems is given in Sections 3 -- 7.

\subsection{Functors $\CC$ and $\CL$}
\label{CLfirst}

Let us study some basic properties of the adjoint functors $\CL$ and $\CC$
defined above. We start with the following Lemma whose proof can be
found in~\cite{haha} (note that the general claim of~\cite{hla}, 3.6.12 
contains an error).
\subsubsection{}
\begin{lem}{U=U}(see~\cite{haha}, 6.8.5) Let $\fg$ be a dg Lie algebra over a 
commutative ring $k\supseteq\Bbb{Q}$.
It can be obviously considered as a strong homotopy Lie algebra, i.e., 
a $\CS$-algebra where $\CS$ is the standard Lie operad (see~\cite{hla} and 
also~\ref{operads}). Suppose that $\fg$ is $k$-flat. 
Then the natural map
$$ U(\CS,\fg)\to U(\Lie,\fg)$$
is a quasi-isomorphism.
\end{lem}
\begin{pf} Recall shortly the reasoning. The functor $U(\CS,\_)$ carries 
quasi-isomorphisms of flat dg Lie algebras into quasi-isomorphisms since 
$\CS$ is a cofibrant operad.

The functor $U(\Lie,\_)$ carries quasi-isomorphisms of flat dg Lie algebras
into quasi-isomorphisms by PBW (see~\cite{rht}, App. B). Finally, the map
$U(\CS,\fg)\to U(\Lie,\fg)$ is a quasi-isomorphism for a cofibrant $\fg$
by the Comparison theorem~\cite{haha}, 5.5.1.

Taking into account that cofibrant dg Lie algebras are flat, we get the result.
\end{pf}

The following \Prop{CL} has a very important filtered analog
 --- see~\ref{CL-filt} below.

\subsubsection{}
\begin{prop}{CL}

(1) The functor $\CC$ preserves quasi-isomorphisms.

(2) The adjunction maps $i_X:X\to\CC\CL(X)$ and 
$p_{\fg}:\CL\CC(\fg)\to \fg$ are quasi-isomorphisms.

(3) The restriction of $\CL$ to the subcategory $\dgcu^{\geq 0}(k)$
of non-negatively graded coalgebras preserves quasi-isomorphisms.
\end{prop}
\begin{pf}

Step 1. Let us prove first that the map $p_{\fg}:\CL\CC(\fg)\to \fg$ is a 
quasi-isomorphism.
Consider the Lie algebra $\fg$ as an algebra over the operad $\CS$ 
as above. According
to~\Lem{U=U}, the natural map of the enveloping algebras
$U(\CS,\fg)\to U(\fg)$ is a quasi-isomorphism. Now, one has
an isomorphism $U(\CS,\fg)=U(\CL\CC(\fg))$ and then the 
adjunction map  $p_{\fg}$ is quasi-isomorphism by the PBW theorem~\cite{rht}, 
App. B.
 
Step 2. Now we check that $\CC$ preserves quasi-isomorphisms. 
The coalgebra $\CC(\fg)$ admits an increasing filtration natural in $\fg$ 
so that the associated graded pieces are $S^n(\fg[1])$. This clearly 
implies the claim.

Step 3. Now we can prove that the map  $i_X:X\to\CC\CL(X)$ is a 
quasi-isomorphism for any $X\in\dgcu(k)$. In fact, the map $i_X$
admits a natural splitting $q_X:\CC\CL(X)\to X$ as a map of complexes:
$Y=\CC\CL(X)$ as a graded vector space takes form
$$Y=S(F(X[-1])[1])$$
where $S$ is the symmetric algebra and $F$ is the free Lie algebra functor.
This defines the projection $q_X:Y\to F(X[-1])[1]\to X$ which is compatible 
with the differentials and splits $i:X\to Y$. Now consider the diagram
\begin{center}
$$
\begin{CD}
X @>{i_X}>> \CC\CL(X) @>{q_X}>> X         \\
@V{i_X}VV    @VV{\CC\CL(i_X)}V  @VV{i_X}V \\
\CC\CL(X) @>>{i_{\CC\CL(X)}}> \CC\CL\CC\CL(X) 
@>>{q_{\CC\CL(X)}}> \CC\CL(X) \\
\end{CD}
$$
\end{center}

The left square in it is commutative by the general nonsense of adjoint
functors; the right square is commutative since $q_X$ is functorial in $X$.
The map $\CL(i_X)$ is a quasi-isomorphism since it is split by the
quasi-isomorphism $p_{\CL(X)}$ by Step 1. Then, by Step 2,
the map $\CC\CL(i_X)$ is also a quasi-isomorphism, and therefore
its retract $i_X$ is a quasi-isomorphism as well.

Step 4. The claim (3) follows by a standard spectral sequence argument.

\end{pf}

\section{Filtered world and graded world}
\label{filtered}
In this Section we prove a filtered analog of~\Lem{U=U} and of
~\Prop{CL} --- see Proppositions ~\ref{U=U-f},~\ref{CL-filt} below. 
For this we need a number of new categories and functors and a 
well-known Rees trick which allows one to reduce some filtered objects to
graded objects over the polynomial ring --- see~\ref{rees}.

\subsection{Filtered world}
\label{f}
\subsubsection{Definitions}
Here $k$ is a base commutative ring. A filtered $k$-module $V$ is a
collection $V=\{V_i\},\ i\in\Bbb{Z}$ with $V_i\subseteq V_{i+1}$ and 
$V=\cup V_i$. The category of filtered $k$-modules is denoted $\Modf(k)$.

A filtered complex is a complex in $\Modf(k)$. The category of filtered
complexes will be denoted in the sequel $CF(k)$ instead of $C(\Modf(k))$.

The category $\Modf(k)$ admits a tensor structure given by the formula

$$ (X\otimes Y)_n=\sum_{p+q=n}\im(X_p\otimes Y_q\to X\otimes Y).$$
This tensor structure induces a tensor structure on $CF(k)$.

The functor $\#:CF(k)\to C(k)$ forgetting the filtration preserves the tensor
structure.

\subsubsection{}
There is an obvious functor 
$$\tau:C(k)\to CF(k)$$
given by $\tau(X)_{-1}=0;\ \tau(X)_n=X, \text{ for } n\geq 0.$
The functor $\tau$ preserves the tensor structure. Thus $\tau$ induces
a functor $\tau:\Op(C(k))\to\Op(CF(k))$. For an operad $\CO\in\Op(C(k))$
we denote by $\Algf(\CO)$ (instead of $\Alg(\tau(\CO))$)  the category
of filtered $\CO$-algebras. 

\subsubsection{}
\begin{exa}{filtered-dglie-dgc}
We write $\dglf(k)$ instead of $\Algf(\Lie)$ for the category of filtered 
dg Lie algebras. Explicitly, such an algebra is a filtered complex 
$\fg=\{\fg_i\}$ with a Lie bracket satisfying $[\fg_i,\fg_j]\subseteq
\fg_{i+j}$. Similarly, we write $\dgcf$ for the category of filtered
cocommutative coalgebras. Its objects are filtered complexes $X=\{X_i\}$
endowed with a cocommutative comultiplication satisfying
$$ \Delta(X_n)\subseteq\sum_{p+q=n}X_p\otimes X_q.$$
\end{exa}

\subsubsection{}
\label{env-f}
Fix $\CO\in\Op(C(k))$ and let $A\in\Algf(\CO)$. The filtration on $A$
induces a natural filtration on the enveloping algebra $U(\CO,A^{\#})$
defined as follows. Recall that $U(\CO,A^{\#})$ is a quotient of the 
``$\CO$-tensor algebra''
$$ T(\CO,A^{\#})=\bigoplus_{n\geq 0}\CO(n+1)
\otimes_{\Sigma_n}(A^{\#})^{\otimes n}.$$
We endow $T(\CO,A^{\#})$ with the tensor product filtration and $U(\CO,A^{\#})$
with the quotient filtration. One easily sees that this defines a
filtered associative dg algebra which will be denoted in the sequel by
$U(\CO,A)$.

\subsection{Graded world}
\label{g}
Let $R$ be a graded commutative (not super-commutative!) $k$-algebra
and let $R^{\#}$ be the corresponding commutative algebra with forgotten
grading.

Denote by $\Modg(R)$ the category of graded $R$-modules. 
It has an obvious tensor structure with the isomorphism
$$ M\otimes_RN\to N\otimes_RM$$
given by the formula $m\otimes n\mapsto n\otimes n$ (no signs involved).

One has the obvious forgetful
functor $\#:\Modg(R)\to\Mod(R^{\#})$ which also preserves the tensor structure.

Denote $CG(R)=C(\Modg(R))$. The forgetful functor defines a tensor functor
$$\#:CG(R)\to C(R^{\#}).$$

\subsubsection{} Tensoring by $R$ defines a tensor functor
$$\tau:\Mod(k)\to\Modg(R).$$
This allows one, for any operad $\CO\in\Op(C(k))$, consider the category
of $\tau(\CO)$-algebras. This latter will be denoted $\Algg(\CO,R)$ or just
$\Algg(\CO)$ (this will not lead to a confusion).

The enveloping algebra $U(\CO,A)$ of $A\in\Algg(\CO)$ is defined in a standard
way as in~\cite{hla} using the tensor structure on $CG(R)$.

We will need the following graded analog of~\Lem{U=U}.

\subsubsection{}
\begin{prop}{U=U-g}
Let $\fg$ be a flat graded dg Lie algebra over $R$. 
Then the natural map 
$$ U(\CS,\fg)\to U(\Lie,\fg)$$
is a graded quasi-isomorphism.
\end{prop}
\begin{pf}Since the forgetful functor $\#:\Modg(R)\to\Mod(R^{\#})$
is exact, and since a map $f:X\to Y$ is a graded quasi-isomorphism iff $f^{\#}$
is a quasi-isomorphism, the result immediately follows from~\Lem{U=U}.
\end{pf}

\subsection{Rees functor}
\label{rees}
$$ $$

{\em From now on $k$ is a field of characteristic zero and  $R=k[t]$ with 
$\deg(t)=1$.
}

The Rees functor
$$\rho:\Modf(k)\lra\Modg(R)$$
is defined by the formula
$$\rho(V)=\sum V_it^i\subseteq \tau(V)=V\otimes R.$$

\subsubsection{}
\begin{lem}{}
1. Rees functor preserves the tensor structure.

2. One has $\rho\tau=\tau$ (two different $\tau$, from ~\ref{f} and 
from~\ref{g}, are involved).
\end{lem}
\begin{pf}Straighforward.
\end{pf}

\subsubsection{}
\begin{cor}{}The Rees functor induces a functor
$$\rho:\Algf(\CO)\to\Algg(\CO).$$
\end{cor}

\subsubsection{}
The Rees functor $\rho$ identifies the category $\Modf(k)$ with the full
subcategory of $\Modg(R)$ consisting of graded torsion-free (=flat) 
$R$-modules. The functor $\rho$ admits a left adjoint functor 
$\phi:\Modg(R)\to\Modf(k)$ defined by the formulas
$$\phi(M)=\dirlim M_n=M/(1-t)M;\ \phi(M)_n=\im(M_n\to\phi(M)).$$

\subsubsection{}
\begin{prop}{u=fur}
Let $\CO\in\Op(C(k)),\ A\in\Algf(\CO)$. The filtered enveloping algebra
can be calculated by the formula
$$ U(\CO,A)=\phi(U(\CO,\rho(A))).$$
\end{prop}
\begin{pf}
The total space of $\phi(U(\CO,\rho(A)))$ is equal to 
$$U(\CO,\rho(A))\otimes_RR/(1-t)R)=U(\CO,\rho(A)\otimes_RR/(1-t)R=
U(\CO,A).$$
To identify the filtration, recall that $U(\CO,\rho(A))_n$ is the image
of the $n$-th component of the tensor algebra $T(\CO,\rho(A))$ which is
an image of
$$ \bigoplus_{i_1+\ldots+i_k=n} \CO(k+1)\otimes A_{i_1}\otimes\ldots
\otimes A_{i_k}.$$

But this coincides with the definition of the filtration on $U(\CO,A)$ 
as in~\ref{env-f}.

\end{pf}

\subsubsection{}
\begin{cor}{ru=ur}
Let $\CO=\Lie\text{ or }\CS$. Then for any $A\in\Algf(\CO)$ one has
$$ \rho(U(\CO,A))=U(\CO,\rho(A)).$$
\end{cor}
\begin{pf}
Having in mind~\Prop{u=fur}, it is enough to check that $U(\CO,\rho(A))$ is 
torsion-free for $\CO=\Lie\text{ or }\CS$. 

If $\CO=\Lie$ the claim follows from the PBW theorem.
In the second case  $\CO=\CS$ is semi-free, so that the corresponding 
enveloping algebra is a tensor algebra which has no torsion.
\end{pf}

\subsubsection{} Note the following nice (surely well-known) generalization
of the PBW. 

\begin{cor}{genpbw}
Let $\fg$ be a filtered Lie algebra. Then the associated graded of the
filtered enveloping algebra $U(\fg)$ is isomorphic to the enveloping
algebra of the associated graded Lie algebra.
\end{cor}
\begin{pf}The passage to associated graded module is the composition of the
Rees functor with the base change $R\to R/(t)$.
\end{pf}

\subsubsection{}
Comparing~\Cor{ru=ur} with~\Prop{U=U-g} we get the following filtered version
of~\ref{U=U}.

\begin{prop}{U=U-f}
Let $\fg$ be a filtered dg Lie algebra over $k$. The natural map
$$ U(\CS,\fg)\to U(\Lie,\fg)$$
is a filtered quasi-isomorphism.
\end{prop}

\subsubsection{}
Note also the following filtered version of PBW.

\begin{lem}{PBW-f}
Let $\fg$ be a filtered dg Lie algebra over $k$. The symmetrization map
$S(\fg)\to\U(\Lie,\fg)$ is a filtered isomorphism. 
\end{lem}
\begin{pf}Use the usual PBW for the dg $R$-Lie algebra $\rho(\fg)$.
\end{pf}

\subsection{A filtered version of~\Prop{CL}}

Let $\fg$ be a filtered dg Lie algebra. The coalgebra $\CC(\fg)$ endowed
with the induced filtration is a filtered unital coalgebra.
{\em Note that filtered unital coalgebras admit two filtrations, the one
being the given filtration, and the second being defined by the unit.}
In the same way the functor $\CL$ sends filtered unital coalgebras to 
filtered Lie algebras.

\subsubsection{}
\begin{defn}{admissible}1. A unital filtered coalgebra 
$X=\{X_i\}\in\dgcf(k)$ is called {\em admissible}
(or, in other words, $\{X_i\}$ is an admissible filtration on $X$)
if $X_{-1}=0$ and $X_0=k\cdot 1$.

2. A filtered dg Lie algebra $\fg=\{\fg_i\}\in\dglf(k)$ is admissible if
$\fg_0=0$.
\end{defn}

\subsubsection{}
\begin{note}{}
If $X$ is an admissible coalgebra, one has
$$ \Delta(x)-1\otimes x-x\otimes 1\in X_n$$
whenever $x\in X_{n+1}$. This follows from the formula 
$$ (1\otimes\epsilon)\Delta=(\epsilon\otimes 1)\Delta=\id.$$
\end{note}

The category of admissible filtered coalgebras is denoted by $\dgca(k)$,
and the category of admissible dg Lie algebras is $\dgla(k)$.

\subsubsection{}
\begin{prop}{CL-filt}1. The functors $\CL$ and $\CC$ define a pair of
adjoint functors
$$ \CL:\dgca(k)\rlarrows\dgla(k).$$
2. The functor $\CC$ preserves filtered quasi-isomorphisms. 

3. The adjunction maps $i_X:X\to\CC\CL(X)$ 
and $p_{\fg}:\CL\CC(\fg)\to\fg$ are  filtered quasi-isomorphisms.
\end{prop}
\begin{pf}
For $X\in\dgca(k)$ and $\fg\in\dgla(k)$ the sets
$\Hom(\CL(X),\fg)$ and $\Hom(X,\CC(\fg))$ coincide with the collection
of filtration preserving twisted cocycles from $X$ to $\fg$. This proves
the first claim.

Now  we can proceed as in the proof of~\ref{CL}.

{\em Step 1.} Let $\fg\in\dgla(k)$. To check that the adjunction map
$p_{\fg}:\CL\CC(\fg)\to\fg$ is a filtered quasi-isomorphism, note that
$U(\CS,\fg)=U(\Lie,\CL\CC(\fg))$ as filtered algebras so that~\Prop{U=U-f}
together with~\ref{PBW-f} give what we need.

{\em Step 2.} Exactly as in the proof of~\ref{CL}, $\CC$ preserves filtered
quasi-isomorphisms since the coalgebras $\CC(\fg)$ admits an increasing
filtration natural in $\fg$ with the associated graded pieces $S^n(\fg[1])$.

{\em Step 3.} Consider now the diagram of the Step 3 of the proof of~\ref{CL}.
Since all the maps involved preserve the filtrations, and since the map
$p_{\CL(X)}$ is a filtered quasi-isomorphism by Step 1, the map $\CL(i_X)$
and, therefore, its retract $i_X$, are also filtered quasi-isomorphisms.

This proves the Proposition.
\end{pf}

Filtered quasi-isomorphisms of admissible coalgebras are useful because
of the following

\subsubsection{}
\begin{prop}{when-we}
Let $f:X\to Y$ be a filtered quasi-isomorphism of admissible 
coalgebras. Then $\CL(f)$ is a quasi-isomorphism.
\end{prop}
\begin{pf}Since the functor $\CL$ commutes with colimits and passing to
cohomology commutes with filtered colimits, the claim can be proven by 
induction. Suppose, by the inductive hypothesis, that the map
$\CL(f_n):\CL(X_n)\to\CL(Y_n)$ is a quasi-isomorphism. Denote $M=X_{n+1}/X_n$
and $N=Y_{n+1}/Y_n$. Choosing compatible splittings for the short exact 
sequences $X_n\to X_{n+1}\to M$ and
$Y_n\to Y_{n+1}\to N$, we get compatible  maps $\alpha:M[-2]\to\CL(X_n)$
and $\beta:N[-2]\to\CL(Y_n)$. Since $M$ and $N$ are quasi-isomorphic,
and $\CL(X_{n+1})=\CL(X_n)\langle M,\alpha\rangle$,
$\CL(Y_{n+1})=\CL(Y_n)\langle N,\beta\rangle$ 
(notation of~\cite{haha}, Sect.~1), Proposition follows from the
following  lemma in the spirit of~\cite{haha}, 5.3.3.
\end{pf}

\subsubsection{}
\begin{lem}{inv-of-ext}
Let $\CO$ be an operad in $C(k)$. Suppose a commutative square
$$
\begin{CD}
$$ M@>\alpha>>X \\
@VgVV      @VfVV \\
   N@>\beta>>Y \\
\end{CD}
$$
is given so that $f:X\to Y$ is a quasi-isomorphism of cofibrant $\CO$-algebras
and $g:M\to N$ is a quasi-isomorphism of complexes. Then the induced
map of $\CO$-algebras $X\langle M,\alpha\rangle\to Y\langle N,\beta\rangle$
is a quasi-isomorphism.
\end{lem}
\begin{pf}
One easily reduces the claim to the case $M=N=k\cdot e$ is generated
by an only element $e$. Since $X$ and $Y$ are cofibrant, $f$ is homotopy 
equivalence. This means that there exists a map $g:Y\to X$ homotopically 
inverse to $f$. Let $x=\alpha(e), y=\beta(e)=f(x)$. The difference $x-g(y)$
is obviously a boundary in the complex $X$; write it as $x-g(y)=du,\ u\in X$.

Choose a path diagram $X\overset{i}{\to} X^I\rra X$ so that $X^I$ is a
standard acyclic cofibration (i.e. has form $X\coprod F(V)$ with a 
contractible $V$) and a map $\Phi:X\to X^I$ which realizes a homotopy
between $\id_X$ and $gf:X\to X$. The map 
$$i': X\langle e; de=x\rangle\to X^I\langle e; de=i(x)\rangle$$
is a quasi-isomorphism since $X^I=X\coprod F(V)$. Since $i(x)$ and $\Phi(x)$
represent one and the same cohomology class in $X^I$, there is an obvious
isomorphism between $X^I\langle e; de=i(x)\rangle$ and 
$X^I\langle e; de=\Phi(x)\rangle$. This implies that the arrows $p'_i$
in the diagram below are quasi-isomorphisms. 
\begin{center}
\begin{picture}(13,3)
   \put(0,1){\makebox(3,1){$X\langle e; de=x\rangle$}}
   \put(5,1){\makebox(3,1){$X^I\langle e; de=\Phi(x)\rangle$}}
   \put(10,0){\makebox(3,1){$X\langle e; de=x\rangle$}}
   \put(10,2){\makebox(3,1){$X\langle e; de=x+du\rangle$}}
   \put(3.2,1.5){\vector(1,0){1.5}}
   \put(8.4,1.8){\vector(2,1){1}}
   \put(8.4,1.2){\vector(2,-1){1}}

   \put(3.2,1.5){\makebox(1.5,1){$\scriptstyle \Phi'$}}
   \put(8.4,2.1){\makebox(.5,.5){$\scriptstyle p'_1$}}
   \put(8.6,.3){\makebox(.5,.5){$\scriptstyle p'_2$}}

\end{picture}
\end{center}
Thus, the map $\Phi'$,
and therefore  the composition
$$p'_1\Phi': X\langle e; de=x\rangle\to X\langle e; de=x+du\rangle$$
induced by the map $gf:X\to X$ is a quasi-isomorphism. The same is of course
true for the morphism $fg:Y\to Y$. This implies that the map 
$$f':X\langle e; e=x\rangle\to Y\langle e; e=y\rangle$$
is a quasi-isomorphism.
\end{pf}

\section{Proofs of the theorems}

\subsection{}

The proof of~\Thm{cmc-uc} is based on the following Key Lemma whose
proof we postpone until Section~\ref{proof-key-lemma}.

\subsubsection{}
\begin{klem}{key-lemma}
Given $X\in\dgcu(k)$, let $f:\fg\to\CL(X)$ be a surjective map in $\dgl(k)$.
Consider the cartesian diagram
$$
\begin{CD}
Z@>j>> \CC(\fg) \\
@VVV    @VV{\CC(f)}V   \\
X@>{i_X}>>\CC\CL(X)\\
\end{CD}
$$
in $\dgcu(k)$. Then $\CL(j):\CL(Z)\to\CL\CC(\fg)$ is a quasi-isomorphism.
\end{klem}

The lemma is very close to~\cite{rht}, II.5.6. Its proof in our context
uses the technicalities of Section~\ref{filtered}.

\subsection{Proof of~\Thm{cmc-uc}}
To prove \Thm{cmc-uc} we use \Lem{key-lemma} and follow the proof 
of Theorem II.5.2 of~\cite{rht}.

\subsubsection{Limits and colimits in $\dgcu(k)$} 
\label{lcl-dgcu}
The functor $\#:\dgcu(k)\to C(k)$
defined by the formula $\#(X)=\ol{X}$, commutes with colimits. This gives a 
obvious construction of arbitrary colimits in $\dgcu(k)$.

Finite products in $\dgcu(k)$ correspond to tensor products of the underlying
complexes; also kernels of a pair of maps in $\dgcu(k)$ are the same as
in $C(k)$. This proves the property CM1 ---  see~\cite{rht}, II.1.

\subsubsection{Cofibrations in $\dgcu(k)$}

$ $

\begin{lem}{cof}
Let $f:X\to Y$ be a cofibration (resp., an acyclic cofibration) in $\dgcu(k)$.
Then $\CL(f)$ is a cofibration (resp., an acyclic cofibration) in $\dgl(k)$.
\end{lem}
\begin{pf}
Let  $f:X\to Y$ be injective, $\{Y_i\}$ be the canonical filtration of $Y$,
$Z_i=f(X)+Y_i\subseteq Y$ are subcoalgebras in $Y$. Since $Y=\dirlim Z_i$
and $\CL$ commutes with colimits, it is enough to check that 
$\CL(Z_i)\to\CL(Z_{i+1})$ is a cofibration. Actually, since $Z_i/Z_{i+1}$
is primitive, the Lie algebra $\CL(Z_{i+1})$ is obtained from $\CL(Z_i)$
by ``joining variables to kill cycles'' procedure. 

The claim about acyclic cofibrations follows from the above and from
the definition of weak equivalences in $\dgcu(k)$.
\end{pf}

\subsubsection{ Fibrations in $\dgcu(k)$}

$ $

\begin{lem}{fib}
Let $f:\fg\to\fh$ be a surjective map (resp., a surjective quasi-isomorphism)
in $\dgl(k)$. Then $\CC(f)$ is a fibration (resp., an acyclic fibration)
in $\dgcu(k)$.
\end{lem}
\begin{pf}
If $f$ is surjective, $\CC(f)$ is a fibration by~\Lem{cof} and 
the adjointness of $\CL$ and $\CC$. If, moreover, $f$ is a quasi-isomorphism,
$\CC(f)$ is a weak equivalence by~\Prop{CL}(2).
\end{pf}

\subsubsection{}
The properties (CM2), (CM3) are obvious. Also the lifting  property (CM4)(ii)
is valid by definition of fibrations in $\dgcu(k)$.

(CM5)(i). Given a map $f:X\to Y$ in $\dgcu(k)$ let $\CL(f)=p i$
be a decomposition of $\CL(f)$ with a cofibration $i:\CL(X)\to\fg$ and an 
acyclic fibration $p:\fg\to\CL(Y)$.

Let $Z=Y\times_{\CC\CL(Y)}\CC(\fg)$. According to~\Lem{key-lemma}
and~\Prop{CL}, the map $Z\to Y$ is a weak equivalence. It is also a fibration
since it is obtained by a base-change from the fibration $\CC(p)$ in 
$\dgcu(k)$. Now, the induced map $X\to Z$ being obviously injective,
we get (CM5)(i).

(CM4)(i). If $f:X\to Y$ is an acyclic fibration, the already proven property
(CM5)(i) gives a decomposition $f=q j$ where $j$ is an acyclic 
cofibration and $q$ is obtained by a base-change from a map $\CC(p)$
with $p$ an acyclic fibration in $\dgl(k)$. Adjointness of $\CL$ and $\CC$
immediately gives that $\CC(p)$ has a RLP with respect to all cofibrations.
Therefore $q$ admits the same property. Since $j$ admits a LLP with respect to
$f$, we get that $f$ is a retract of $q$ and therefore, it also satisfies RLP
with respect to all cofibrations.

(CM5)(ii). Let $f:X\to Y$ and let  $\CL(f)=p i$ be a decomposition with 
a fibration $p:\fg\to\CL(Y)$  and an acyclic cofibratiton $i:\CL(X)\to\fg$.
According to \Lem{key-lemma} the map $j:Z\to\CC(\fg)$ is a weak equivalence
and the map $Z\to Y$ is a fibration where 
$Z=Y\times_{\CC\CL(Y)}\CC(\fg)$. This immediately implies that
the map $X\to Z$ is an acyclic cofibration.

Therefore, we proved $\dgcu(k)$ admits a CMC structure. The simplicial 
structure on $\dgcu(k)$ and the proof of the axiom (SM7) will be provided
in Section~\ref{simpl}.
 
\subsection{Proof of \Thm{2}}

Now \Thm{2} follows immediately from the general Theorem II.1.4 
of~\cite{rht} and from \Prop{CL}.

\section{Proof of the Key \Lem{key-lemma}}
\label{proof-key-lemma}

In this Section we prove the Key \Lem{key-lemma}. 

Endow $X$ with the canonical filtration and $\CL(X)$ with the induced 
filtration.  Let $\fa=\Ker(f:\fg\to\CL(X))$.  Define an admissible filtration
on $\fg$ by setting $\fg_n=f^{-1}(\CL(X)_n)$ for $n>0$. This induces
admissible filtrations on $\CC(\fg)$ and on $\CC(\CL(X))$.

According to ~\Prop{CL-filt} $i_X:X\to\CC\CL(X)$
is a filtered quasi-isomorphism. 
Define finally a filtration on $Z$ by the formula
$Z_n=j^{-1}(\CC(\fg)_n)$. 

According to~\Prop{when-we}
it is enough to check that $j:Z\to\CC(\fg)$ is a filtered quasi-isomorphism.

Let us describe more explicitly the filtrations involved. Forget about the 
differentials. Choose a graded Lie algebra splitting $s:\CL(X)\to\fg$ of 
$f$. This defines isomorphisms (not preserving the differentials)
\begin{equation}
\CC(\fg)\iso\CC\CL(X)\otimes\CC(\fa)
\label{filtC(g)}
\end{equation}
 and
\begin{equation}
Z\iso X\otimes\CC(\fa)
\label{filtZ}
\end{equation}
 of filtered coalgebras, the filtration on $ \CC(\fa)$
being the standard one. In fact, the first isomorphism obviously preserves
filtrations, and the
second one preserves the filtrations because of the equality 
$X_n=i_X^{-1}(\CC\CL(X)_n)$.  

The isomorphism~(\ref{filtZ}) can be rewritten as
\begin{equation}
Z_n\iso\bigoplus_r X_{n-r}\otimes S^r(\fa[1])
\label{filtZ1}
\end{equation}
and similarly 
\begin{equation}
\CC(\fg)_n\iso\bigoplus_r\CC\CL(X)_{n-r}\otimes S^r(\fa[1]).
\label{filtC(g)1}
\end{equation}
Unfortunately, the isomorphisms~(\ref{filtZ1}),~(\ref{filtC(g)1})
are not compatible with the differentials. To overcome this minor 
difficulty, we define a double filtration on the complexes involved
so that the associated graded guys will be already isomorphic as
complexes. We will write formulas only for the filtration on $Z$ and
on $Z_n$, the formulas for $\CC(\fg)$ being obtained by
substitution of $X_n$ with $\CC\CL(X)_n$. Here are the formulas.

\begin{equation}
F^q_p=\bigoplus_{r\geq q}X_p\otimes S^r(\fa[1])
\label{FZ}
\end{equation}

\begin{equation}
F^q_p(Z_n)=F^q_p\cap Z_n=
\bigoplus_{n\geq r\geq q}X_{\min(p,n-r)}\otimes S^r(\fa[1])
\label{FZn}
\end{equation}

The filtrations are increasing on $p$ and decreasing on $q$. The 
filtration~(\ref{FZn}) is finite. Its $(p,q)$-graded piece vanishes
for $p+q>n$ and is otherwise isomorphic to $X_p/X_{p+1}\otimes S^q(\fa[1])$
as a complex.

Associated graded pieces of the corresponding filtration for $\CC(\fg)$
have form  $\CC\CL(X)_p/\CC\CL(X)_{p+1}\otimes S^q(\fa[1])$. 
 
The $(p,q)$-graded piece of the map $j_n:Z_n\to\CC(\fg)_n$ takes form
$$ X_p/X_{p+1}\otimes S^q(\fa[1])\lra
\CC\CL(X)_p/\CC\CL(X)_{p+1}\otimes S^q(\fa[1])
$$
which obviously a quasi-isomorphism by~\Prop{CL-filt}.

Key Lemma is proven.

\section{Simplicial structure on $\dgcu(k)$}
\label{simpl}

In this Section we define a simplicial structure on the category $\dgcu(k)$
of dg unital coalgebras over a field $k$ of characteristic zero and check
the axiom (SM7) --- see Introduction.

\subsection{Functional spaces for unital coalgebras}
Recall  (cf.~\cite{bg}) that the functor of polynomial differential forms
\begin{equation}
\Omega:\simpl\to\dga(k)
\end{equation}
is the one defined uniquely by its values on the standard simplices
$$\Omega(\Delta^n)=\Omega_n=
k[t_0,\ldots,t_n,dt_0,\ldots,dt_n]/(\sum t_i-1,\sum dt_i)$$
and by the property that $\Omega$ commutes with colimits.

\subsubsection{}
For any commutative dg algebra $\Omega\in\dga(k)$ 
tensoring by $\Omega$ defines a functor 
\begin{equation}
 \Omega\otimes\_:\dgcu(k)\to\dgcu(\Omega).
\label{bch}
\end{equation}

Therefore, the following definition makes sence.

\subsubsection{}
\begin{defn}{func}
Let $X,Y\in\dgcu(k)$. The simplicial set $\uhom(X,Y)$ is defined by
the formula
$$\uhom(X,Y)_n=\Hom_{\dgcu(\Omega_n)}(\Omega_n\otimes X, \Omega_n\otimes Y),$$
the faces and the degeneracies being defined in an obvious way. 
\end{defn}

Note the following

\subsubsection{} 
\begin{lem}{lims-bch}
The functor~(\ref{bch}) commutes with colimits and with finite limits
\end{lem}
\begin{pf}
The claim about colimits is obvious. In order to prove that~(\ref{bch})
commutes with finite limits we check separately the case of a product
of two coalgebras and that of kernel of a couple of maps. This follows from 
the description of limits in~\ref{lcl-dgcu}. 
\end{pf}

\Lem{lims-bch} immediately implies the following

\subsubsection{}
\begin{cor}{good-Hom}1. The functor $\uhom(X,\_):\dgcu(k)\to\simpl$
commutes with finite limits.

2.  The functor $\uhom(\_,Y):\dgcu(k)^{\op}\to\simpl$
carries arbitrary colimits to limits.
\end{cor}

One has the following standard fact.
\subsubsection{}
\begin{lem}{func-repr}(see~\cite{bg},~Lemma 5.2,~\cite{haha}, 4.8.3) There is 
a natural in \\
$S\in\simpl$ morphism
$$\Phi(S):\Hom_{\dgcu(\Omega(S))}(\Omega(S)\otimes X,\Omega(S)\otimes Y)\iso
\Hom(S,\uhom(X,Y))$$
which is a bijection provided $S$ is finite.
\end{lem}
\begin{pf} The proof of~\cite{haha}, 4.8.3 is applicable here.
\end{pf}

\subsubsection{}
\begin{lem}{LadjC-simpl}
The adjoint functors $\CC$ and $\CL$ induce an isomorphism
$$\uhom(X,\CC(\fg))\iso\uhom(\CL(X),\fg)$$
of simplicial sets for every $X\in\dgcu(k),\ \fg\in\dgl(k)$.
\end{lem}
\begin{pf}Repeats the standard argument of~\Thm{LadjC} substituting the base
category $C(k)$ with $\Mod(\Omega_n)$.
\end{pf}

\subsection{Property (SM7)}

\subsubsection{}
\begin{prop}{SM7}
 Let $i:A\to B$ be a cofibration and $p:X\to Y$ be a fibration in $\dgcu(k)$.
Then the map of simplicial sets
\begin{equation}
 \pi(i,p):\uhom(B,X)\to\uhom(A,X)\times_{\uhom(A,Y)} \uhom(B,Y)
\label{map7}
\end{equation}
is a Kan fibration. If, moreover, either $i$ or $p$ is a weak equivalence,
then $\pi(i,p)$ is an acyclic Kan fibration.
\end{prop}

\subsubsection{}
In what follows we will say that a pair of maps $(i:A\to B,p:X\to Y)$
satisfies (SM7) if the map $\pi(i,p)$ from~(\ref{map7}) is a Kan fibration,
acyclic if one of $(i,p)$ is a weak equivalence. To prove Proposition, we will
show step by step that any pair  $(i:A\to B,p:X\to Y)$ such that
$i$ is a cofibration and $p$ is a fibration, satisfies (SM7).

{ \em Step 1.} Suppose that $p=\CC(f)$ where $f:\fg\to\fh$ is a surjective map
of dg Lie algebras. Then any pair $(i,p)$ satisfies (SM7) by~
\ref{LadjC-simpl}, \ref{cof} and \cite{haha}, 4.8.4. 

{\em Step 2.} Suppose that a pair $(i,p)$ satisfies (SM7) and let a map
$q:Z\to T$ be obtained from $p:X\to Y$ by a base change $a:T\to Y$. 

Using~\Cor{good-Hom}, we easily see that $\pi(i,q)$ is obtained by a base
change from $\pi(i,p)$.  Therefore the pair $(i,q)$ also satisfies (SM7).

{\em Step 3.} Suppose that a pair $(i,p)$ satisfies (SM7) and let a map
$q:Z\to T$ be a retract of $p:X\to Y$.  Then the map $\pi(i,q)$ 
is a retract of $\pi(i,p)$ and therefore, the pair $(i,q)$ also satisfies
(SM7).

Now Proposition follows from the following lemma.

\subsubsection{}
\begin{lem}{descr-f-dgcu}1. Any fibration in $\dgcu(k)$ can be obtained,
using the operations of retraction and base change, from a map $\CL(f)$
where $f$ is a surjective map of dg Lie algebras.

2. Any acyclic fibration in $\dgcu(k)$ can be obtained,
using the operations of retraction and base change, from a map $\CL(f)$
where $f$ is a surjective quasi-isomorphism of dg Lie algebras.
\end{lem}
\begin{pf}1. Let $f:X\to Y$ be a fibration. Using the maps $i:Y\to\CC\CL(Y)$
and $\CC\CL(f):\CC\CL(X)\to\CC\CL(Y)$, define $Z=Y\times_{\CC\CL(Y)}\CC\CL(X)$
and let $j:Z\to \CC\CL(X),k:X\to Z$  and $g:Z\to Y$ be the obviously defined 
maps. According to the Key Lemma, $j$ is an acyclic cofibration, and therefore,
$k$ is an acyclic cofibration as well. Since $f$ is a fibration, the map $k$ 
splits over $g$ and this gives a presentation of $f$ as a retract of $g$
which is obtained by a base change frm $\CC(\CL(f))$. 

2. If, moreover, $f$ is an acyclic fibration, then $\CC(\CL(f))$ is a
surjective quasi-isomorphism, and it is nothing to prove.

\end{pf}

\section{The nerve of a dg Lie algebra}

\subsection{The nerve and its properties}
\label{nerve.1}
Let $X\in\dgcu(k)$. Choose a fibrant resolution $X\to F$ and define a
functor 
$$ \widetilde{X}:\dgar(k)\to\simpl$$
by the formula
$$ \widetilde{X}(A)=\uhom(A^*,F)$$
where $A^*$ is the unital coalgebra with the unit $A\to k$.
The resulting functor $\widetilde{X}$ does not depend, up to a homotopy, on 
the choice of the resolution. One can get a specific representative for
$ \widetilde{X}$ as follows.

Let  $\fg=\CL(X)$ be the tangent Lie algebra of $X$. Choose
$\CC(\fg)$ to be a fibrant resolution of $X$. This allows one to easily
express the functor $\widetilde{X}$ through the tangent Lie algebra $\fg$.

\subsubsection{}
\begin{defn}{def.nerve}
Let $\fg\in\dgl(k)$. The nerve of $\fg$ is the functor
$$ \Sigma_{\fg}:\dgar(k)\to\simpl$$
defined by the formula
$$\Sigma_{\fg}((A,\fm))_n=\MC(\fm\otimes\Omega_n\otimes\fg).$$
\end{defn}

One has immediately the following
\subsubsection{}
\begin{prop}{nerveis}
For $X\in\dgcu(k)$ the functor $\widetilde{X}:\dgar(k)\to\simpl$
is homotopy equivalent to the nerve $\Sigma_{\CL(X)}$.
\end{prop}
\begin{pf}
According~\Lem{LadjC-simpl}, one has for any $\fg\in\dgl(k),\ A\in\dgar(k)$
\begin{equation}
\Sigma_{\fg}(A)_n=\MC(\Omega_n\otimes\fm\otimes\fg)=
\Hom_{\dgcu(\Omega_n)}(\Omega_n\otimes A^*,\Omega_n\otimes\CC(\fg))= 
\uhom_n(A^*,\CC(\fg)).
\end{equation}
\end{pf}

\subsubsection{}
\begin{prop}{properties}
1. A quasi-isomorphism $f:\fg\to\fh$ of dg Lie algebras
induces a  homotopy equivalence for every $A\in\dgar(k)$
$$\Sigma_f:\Sigma_{\fg}(A)\to\Sigma_{\fh}(A).$$

2. For each $\fg\in\dgl(k)$ the functor $\Sigma_{\fg}$ carries 
quasi-isomorphisms in $\dgar(k)$ to homotopy equivalences. 

3. For each $\fg\in\dgl(k)$ the functor $\Sigma_{\fg}$ carries
surjective maps to Kan fibration. In particular, 
$\Sigma_{\fg}(A)$ is always a Kan simplicial set.

4. $\Sigma_{\fg}$ commutes with finite projective limits.
\end{prop}
\begin{pf}
The claims 1,3,4 follow from~\Thm{} and~\Lem{LadjC-simpl}. 
By~\Prop{CL}, (3), a quasi-isomorphism $f:A\to B$ in $\dgar(k)$ defines
a weak equivalence $f^*:B^*\to A^*$ in $\dgcu(k)$. Then the induced
map $\Sigma_{\fg}(f)$ is a weak equivalence of Kan simplicial sets, hence
a homotopy equivalence. 
\end{pf}

\subsection{Some calculations}
\label{somecalc}
Here we provide some explicit calculations
which help to better understand how the nerve of a dg Lie algebra looks like.
We will use them below in~\ref{moduli-G}.

In this subsection $\fg$ is a nilpotent dg Lie algebra (it substitutes
$\fm\otimes\fg$ from~\ref{def.nerve}), and we denote by $\Sigma(\fg)$
the simplicial set 
$$\Sigma(\fg)_n=\MC(\fg_{\bullet})$$
where the simplicial dg Lie algebra $\fg_{\bullet}=\{\fg_{(n)}\}$ is defined
by the formula $\fg_{(n)}=\Omega_n\otimes\fg$.

\subsubsection{Deligne groupoid}
\label{del-gr}
Recall (cf.~\cite{gm1}) that for a nilpotent dg Lie algebra $\fg$ one 
defines {\em Deligne groupoid} $\Gamma(\fg)$ as follows.

The Lie algebra $\fg^0$ acts on $\MC(\fg)$ by vector fields:
$$ \rho(y)(x)=dy+[x,y]\text{ for }y\in\fg^0, x\in\fg^1.$$
This defines the action of the nilpotent group $G=\exp(\fg^0)$ on the
set $\MC(\fg)$. Then the groupoid $\Gamma(\fg)$ is defined by the formulas
$$ \Ob\Gamma=\MC(\fg)$$
$$ \Hom_{\Gamma}(x,x')=\{g\in G|x'=g(x)\}.$$

\subsubsection{Maurer-Cartan elements of $\fg_{(1)}$}
Let us explicitly describe the set 
$\MC(\fg_{(1)})$. Since $\fg_{(1)}=k[t,dt]\otimes\fg$ we will iterate the
calculation to get the description of $\MC(\fg_{(n)})$.

Write an element $z\in\fg_{(1)}^1$ in the form
$$ z=x+dt\cdot y$$
with $x\in\fg^1[t],\ y\in\fg^0[t]$. Then the Maurer-Cartan equation
is easily seen to be equivalent to the differential equation

$$ \frac{dx(t)}{dt}=dy(t)+[x(t),y(t)]$$
$$  x(0)=x_0$$
where $x_0$ is an element of $\MC(\fg)$.

An element $y\in\fg[t]$ defines a unique polynomial path $g(t)$ in the Lie 
group $G=\exp(\fg^0)$ satisfying the differential equation 
$\dot{g}(t)=g(t)(y(t))$ with the initial condition $g(0)=1$.

Let $\fk$ be the Lie subalgebra $t\fg^0[t]\subseteq\fg_{(1)}^0$ and let
$K=\exp(\fk)$. The above consideration proves the following
\subsubsection{}
\begin{lem}{path}
An element $x$ of $\MC(\fg_{(1)})$ can be uniquely represented in the form
$$ x=g(x_0)$$
where $x_0\in\MC(\fg)\subseteq\MC(\fg_{(1)})$, 
$g\in K\subseteq\exp(\fg_{(1)}^0)$ and the action is defined as in~\ref{del-gr}
--- for the nilpotent Lie algebra $\fg_{(1)}$.
\end{lem}

\subsubsection{}
\label{mc-of-g_n}
Iteratively using~\Lem{path} we can describe the set of Maurer-Cartan elements
of $\fg_{(n)}$ as follows.

Let $\fk_i=t_i\fg_{(i-1)}^0[t_i]$ for $i>0$ be the Lie subalgebra of 
$\fg_{(n)}^0$.
Denote $K_i=\exp(\fk_i)$. These are subgroups of $\exp(\fg_{(n)}^0)$.

\begin{lem}{normal}
Let $0\leq i\leq j\leq n$. Then $K_i$ normalizes $K_j$.
\end{lem}
\begin{pf}
This immediately follows from the inclusion
$$ [\fk_i,\fk_j]\subseteq \fk_j.$$
\end{pf}

Define $G_n=\exp(\fg_{(n)}^0)$ and let $T_n=K_n\cdot K_{n-1}\cdots K_1$
be the subgroup in $G_n$. The lemma above implies that this 
group is the exponent of the Lie algebra $\oplus_{i\geq 1}\fk_i$.
Then one has the following

\subsubsection{}
\begin{prop}{MC-n}
Any element of $\MC(\fg_{(n)})$ can be uniquely presented as $g(x_0)$
where $x_0\in\MC(\fg)$ and $g\in T_n$.
\end{prop}

Note that the simplicial group $G_{\bullet}=\{G_n\}$ acts on the
nerve $\Sigma(\fg)$. \Prop{MC-n} implies that the restriction
\begin{equation}
 G_{\bullet}\times\MC(\fg)\to \Sigma(\fg)
\label{tosigma}
\end{equation}
is surjective. 
One has the following stronger

\subsubsection{}
\begin{cor}{pcs}
The map~(\ref{tosigma})  admits a pseudo-cross section (see~\cite{may}, \S 18).
\end{cor}
\begin{pf} The pseudo-cross section is given as the composition
$$ \Sigma_n(\fg)\iso T_n\times\MC(\fg)\to G_n\times\MC(\fg).$$
This map obviously commutes with the faces $d_i$ for $i>0$ and with
all degeneracies.
\end{pf}

\subsubsection{}
\begin{note}{}
Using the explicit description of $\Sigma(\fg)$ above, one can easily get
the property $(3)$ of~\ref{properties} independently on 
Theorems~\ref{cmc-uc},~\ref{2}.
\end{note}

\section{Remarks and applications}

In this Section we provide some examples, definitions, calculations 
and remarks.

\subsection{Homology of Lie algebras}

Let us give another description of the homology functor
$\#\circ\CC:\dgl(k)\to C(k)$ where the functor
$\#:\dgcu(k)\to C(k)$ is given by the formula $\#(C)=\ol{C}$.

\subsubsection{}
\begin{prop}{c=homology} One has $\#\circ\CC=\Left\Ab$ where 
$\Ab:\dgl(k)\to C(k)$ is defined as
$$ \Ab(\fg)=\fg/[\fg,\fg].$$
\end{prop}
\begin{pf}It suffices to check that the map of complexes 
$\ol{\CC(\fg)}\to\fg/[\fg,\fg]$  is a  quasi-isomorphism provided $\fg$
is standard cofibrant. Consider $\CC(\fg)$ as a bicomplex so that the
horizontal differential $d'$ is defined by the Lie  bracket in $\fg$ and the
vertical differential $d''$ is induced by the differential in $\fg$.
Forget for a moment the differential in $\fg$. Then 
\begin{equation}
\fg=F(V)=\oplus F^n(V)
\label{free}
\end{equation}
is a free Lie algebra over a graded vector space $V$. The differential $d'$
preserves the grading which comes from the presentation~(\ref{free}). Since
$H(\CC,d')=V$ and each homogeneous component of $\CC$ is finite, the
proposition follows.
\end{pf}

The above description allows one to construct easily an example of a 
acyclic coalgebra $X$ such that $\CL(X)$ has non-trivial cohomology
(another example is given in Kontsevich's lectures~\cite{ko}).

\subsubsection{}
\begin{exa}{counter}
Let $\fg$ be the cofibrant Lie algebra having generators $e,f,h$ of degree $0$,
$x,y,z$ of degree $-1$ with the differential given by
$$ de=dh=df=0; dx=[h,e]-2e; dy=[h,f]+2f; dz=[e,f]-h.$$
According to~\Prop{c=homology}, $\CC(\fg)$ is acyclic though $\fg$ is not
(one has $H^0(\fg)=\frak{sl}_2(k)$). One can equally set $X=\CC(\fg)$
and get a non-contractible $\CL(X)$. This counter-example means that there are
quasi-isomorphisms in $\dgcu(k)$ which are not weak equivalences.
\end{exa}

\subsection{Infinitesimals}Look at the first infinitesimal deformations 
corresponding to a dg Lie algebra $\fg$.

\subsubsection{Dual numbers} 
For each $n=0,1,\ldots$ define
$$A_n=k[\epsilon; \deg\epsilon=-n]/(\epsilon^2)\in\dgar(k).$$
This is a $k$-vector space object in the category $\dgar(k)$.

Let us calculate the simplicial vector space $\Sigma_{\fg}(A_n)$. 
Its $i$-simplices
are the Maurer-Cartan elements of $\epsilon\cdot\Omega_i\otimes\fg$ which is 
a dg Lie algebra with zero multiplication. Therefore,
$$ \Sigma_{\fg}(A_n)_i=Z^0(\Omega_i\otimes\fg[1+n]).$$

Using the Dold-Puppe equivalence of categories and the fact that the 
cosimplicial complex $\{\Omega_i\}_{i\in\Bbb{N}}$ is homotopy equivalent to 
the cosimplicial complex of cochains $\{C^*(\Delta^i)\}_{i\in\Bbb{N}}$, 
we obtain the following

\subsubsection{}
\begin{claim}{dual.numbers}$\Sigma_{\fg}(A_n)$ is homotopy equivalent to 
the simplicial abelian group corresponding to the complex 
$\tau^{\leq 0}(\fg[1+n])$.
\end{claim} 

Note the following
\subsubsection{}
\begin{cor}{faithful}If a map $f:\fg\to\fh$ in $\dgl(k)$ induces an equivalence
of the nerve functors, then $f$ is itself a quasi-isomorphism.
\end{cor}

\subsection{Formal spaces}
\subsubsection{} 
\begin{defn}{def.f.space}
A formal stack $X\in\dgcu(k)$ is called a formal space if it is weakly
equivalent to a coalgebra $Y\in\dgcu(k)$ satisfying $Y^i=0$ for $i<0$.
\end{defn}

An equivalent condition: $X$ is a formal space if $H^i(\CL(X))=0$
for $i\leq 0$. According to~\ref{dual.numbers}, this property is equivalent
to the one saying that $\Sigma_{\fg}(A_0)$ is discrete. By an obvious
artinian induction we get

\subsubsection{}
\begin{lem}{f.space}$X\in\dgcu(k)$ is a formal space iff for each $A\in\art(k)$
the simplicial set $\widetilde{X}(A)$ is discrete. 
\end{lem}

\subsubsection{}
Now the two ideas mentioned in the Introduction about formal deformations
in characteristic zero can be formulated as follows.

{\em
Any formal deformation problem in characteristic zero can be described
by a representable functor
$$ F:\dgar(k)\to\simpl.$$
}

Classical deformation problems are often not representable, since in the 
classical picture we see only the $\pi_0$ of the genuine deformation functor.

\begin{defn}{def.coarse}Let $X\in\dgcu(k)$. The classical part of 
$\widetilde{X}$
is the functor
$$\widetilde{X}_{cl}:\art(k)\to\Ens$$
defined by $\widetilde{X}_{cl}(A)=\pi_0(\widetilde{X}(A))$.
\end{defn}

Let $\fg=\CL(X)$. Suppose first that $X$ is a formal space. Put 
$Y=H^0(\CC(\fg))$. Then for any $A\in\art(k)$ one has
$$X_{cl}(A)=\pi_0(\widetilde{X}(A))=\widetilde{X}(A)=\Hom(A^*,X)=\Hom(A^*,Y).$$
This means that the classical part $\widetilde{X}_{cl}$ of a formal space 
$X$ is representable by the coalgebra $Y$. For a general $X\in\dgcu(k)$
one should not expect representability of the classical part. However, 
the functor $\widetilde{X}_{cl}$ admits a hull in the sense of~\cite{sc}.
In fact, choose a complement $V$
in $\fg^1$ to the vector subspace $\im(d:\fg^0\to\fg^1)$ and define
a  1-truncation $\fh$ of $\fg$ by the formulas
\begin{equation}
\fh^i=\left\{\begin{array}{lll}
              0& , & i\leq 0\\
              V& , & i=1 \\
              \fg^i& , & i>1.
             \end{array}
      \right.
\label{1-trunc}       
\end{equation}

Put $Y=H^0(\CC(\fh))$ and define $h_Y(A)=\Hom(A^*,Y)$.
\subsubsection{}

\begin{lem}{coarse}The injection $\fh\to\fg$ induces a smooth morphism
of functors $h_Y\to\widetilde{X}_{cl}$ which is isomorphism
on the tangent spaces.
\end{lem}
\begin{pf} This claim essentially belongs to Goldman-Millson~\cite{gm1},
\cite{gm2} (who considered however only the case of formal spaces).
The tangent spaces to $h_Y$ and to $\widetilde{X}_{cl}$ are
isomorphic to $H^1(\fh)$ and to $H^1(\fg)$ respectively. To check the 
smoothness it is enough to prove that for any surjection 
$f:(B,\fn)\to (A,\fm)$ in $\art(k)$ whose kernel $I$ is annihilated by $\fn$,
the map 
$$h_Y(B)\to h_Y(A)\times_{\widetilde{X}_{cl}(A)}\widetilde{X}_{cl}(B)$$
is surjective.

That is, let $x\in\MC(\fm\otimes V)$ and $y\in\MC(\fn\otimes\fg^1)$
have the same image in $\widetilde{X}_{cl}(A)$. Then, first of all, one can
substitute the element $y$ by an equivalent one, so that the images of $x$
and of $y$ in $\MC(\fm\otimes\fg^1)$ coincide. Then the element $y$ belongs
to $\fn\otimes V$, up to an element from $I\otimes\im(d:\fg^0\to\fg^1)$
which can be killed by the action of $I\otimes\fg^0$. After this correction,
the element $y$ already belongs to $\fn\otimes V$ and it automatically 
satisfies the MC equation.
\end{pf}

Therefore, the coalgebra $Y$ (or, its dual complete local algebra) is 
a hull of $\widetilde{X}_{cl}$. 

We would like to have a direct proof of the uniqueness of $Y$. For this it 
would be enough to prove that $\fh$ does not depend, up to a quasi-isomorphism,
on the choice of $1$-truncation. This is of course so if (as in~\cite{gm2}) 
one supposes that $H^0(\fg)=0$. Unfortunately, we doubt this is true in 
general.

However, due to the general claim of~\cite{sc}, the hull $Y$ is unique up to a 
non-canonical isomorphism provided $H^1(\fg)$ is finite-dimensional.

\subsection{Simply connected rational spaces}
\label{rht-example}
Let $S$ be a simply connected rational space. According to~\cite{rht}, 
it has a dg Lie algebra model $\fg$ which satisfies $\fg^i=0$ for $i\geq -1$
(we keep using degree $+1$ differentials). Therefore, $S$ should define
a formal deformation in our general definition. It looks strangely a bit,
since ``usual'' deformations are described by non-negatively graded Lie 
algebras. The classical part of such a deformation is trivial. However,
one can easily calculate the corresponding to $S$ deformation functor ---
in terms of dg Lie algebra models.

\subsubsection{}
\begin{prop}{rht} Let $S$ have a finite $\Bbb{Q}$-type. 
For any $A\in\dgar(\Bbb{Q})$ the simplicial set 
$\Sigma_{\fg}(A)$ is simply connected and rational. Its Lie algebra model is 
given by $\fm\otimes\fg$.
\end{prop}
\begin{pf}Put $\fh=\fm\otimes\fg$. We wish to check that $\fh$ is a Lie algebra
model for the simplicial set $\Sigma(\fh)$. But this is clear: the coalgebra
model of $\fh$ is $\CC(\fh)$, so the dg algebra model of $\fh$ is the 
dual complex $\CC^*(\fh)$ and the corresponding simplicial set is given
according to ~\cite{bg}, Thm. 9.4, by the formula
$$ n\mapsto \Hom(\CC^*(\fh),\Omega_n)=\MC(\Omega_n\otimes\fh)=\Sigma_n(\fh)$$
since $\fh^i$ are finitely dimensional.
\end{pf}

\subsection{Example: Intersection of subschemes}

A typical example of a formal space which is not a formal scheme is given by a 
non-transversal intersection of subschemes. Let $X,Y\subseteq Z$ be closed
subschemes in a noetherian scheme $Z$, $z\in X\cap Y$. We wish to describe 
the intersection of $X$ and
$Y$ in $Z$ near $z$. Let $A,B,C$ be the local rings of $X,Y,Z$ respectively.

These rings (or the corresponding dual coalgebras $A^*,B^*,C^*$) represent 
functors
$$ F_A,F_B,F_C:\dgar(k)\to \simpl$$
where $k=k(z)$ is the base field.

The functors $ F_A,F_B,F_C$ being defined up to homotopy equivalence, the best
thing we can do is to consider their homotopy fibre product. 

Thus, define the homotopy intersection of $X$ and $Y$ in $z\in Z$ to be the
homotopy fibre product functor $F$ of $F_A$ and $F_B$ over $F_C$.
In order to calculate it, one has to substitute a map $F_A\to F_C$
(or the other one) with a fibration and take the usual fibre product.

For this it suffices to take a cofibrant resolution $\widetilde{A}$ for the 
$C$-algebra $A$ and substitute $F_A$ with $F_{\widetilde{A}}$. The result
will be given by the dg algebra $\widetilde{A}\otimes_CB$ defined canonically
in the corresponding homotopy category, concentrated in the
nonpositive degrees. Its cohomology is given by the formula
$$ H^i(\widetilde{A}\otimes_CB)=\Tor_{-i}^C(A,B)$$
 --- exactly as one could expect.

\subsubsection{}
\begin{rem}{}It is unclear how to define a global object corresponding, say,
to the intersection of two subschemes in a non-affine scheme.
One should probably use a technique suggested by Hirschowitz-Simpson 
in~\cite{hs}.
\end{rem}

\subsection{Example: Quotient by a group action}

Let an algebraic group $G$ over $k$ acts on a dg Lie algebra $\fh$.
Then $G$ acts on each simplicial set $\Sigma_{\fh}(A)$. Define the functor
$$ F:\dgar(k)\to\simpl$$
as the homotopy quotient $F(A)=\Sigma_{\fh}(A)/\widehat{G}_1(A)$
where $\widehat{G}_1$ is the formal completion of $G$ at $1$ so that
$\widehat{G}_1(A)=\exp(\fm\otimes\fg)$.

\subsubsection{}
\begin{prop}{quotient}
The functor $F$ is homotopy equivalent to the nerve
of the semidirect product $\fg\ltimes\fh$.
\end{prop}
\begin{pf}According to~\Prop{properties}, (3), (4), one has a fibration 
$$ f:\Sigma_{\fg\ltimes\fh}(A)\to\Sigma_{\fg}(A)$$
with fibre $\Sigma_{\fh}(A)$. On the other hand, the map~(\ref{tosigma})
gives in our case a fibration
$$\pi: G_{\bullet}(A)\to\Sigma_{\fg}(A)$$
with fibre $\widehat{G}_1(A)$ and contractible total space. Then the
cartesian diagram
\begin{center}
$$
\begin{CD}
 X@>>{\widehat{G}_1(A)}>\Sigma_{\fg\ltimes\fh}(A)  \\
@V{\Sigma_{\fh}(A)}VV                 @VV{\Sigma_{\fh}(A)}V \\
G_{\bullet}(A)@>>{\widehat{G}_1(A)}>\Sigma_{\fg}(A)
\end{CD}
$$
\end{center}
(all the arrows being fibrations marked by the corresponding fibres) 
presents $\Sigma_{\fg\ltimes\fh}(A)$ as a homotopy quotient of  
$\Sigma_{\fh}(A)$ modulo $\widehat{G}_1(A)$.
\end{pf}

\subsection{Example: Moduli of $G$-torsors}
\label{moduli-G}
The last example of a formal stack will be that of moduli of principal 
$G$-bundles.

Let $G$ be an algebraic group over a field $k$ of characteristic zero,
$S$ be a scheme over $k$, $P$ be a $S$-torsor under $G$. We wish
to study formal deformations of $P$. For this we have to define a 
deformation functor
$$F_P:\dgar(k)\to\simpl$$
naturally generalizing the standard functor of formal deformations
$\art(k)\to\Grp$. 

We will proceed as follows. First of all, we define in~\ref{affine} torsors 
with a (affine) dg base. Then (in~\ref{aff.def}) we  construct a simplicial 
set describing the formal deformations of a torsor in the affine case. 

Deformations of a torsor over a non-affine scheme are defined 
in~\ref{non-aff} as an 
appropriate homotopy limit of the deformations over affine bases. 
After this  we explicitly calculate in \ref{calc.triv}--~\ref{triv-d} 
the deformation functor
describing deformations of a trivial torsor over an affine base. 
Finally, using the main result of~\cite{ddg},  we get the final 
answer --- \Thm{torsors}.

Fix some notations. Recall that $\dga(k)=\Alg(\Com(k))$ is the category of 
commutative dg algebras over $k$.

Let $R$ be the Hopf algebra of regular functions on $G$. Then
$\fg=\Der(R,R)$ is the Lie algebra of $G$.

\subsubsection{}
\label{affine}
\begin{defn}{tors}
Let $B\in\dga(k)$. A $B$-torsor under $G$ is a morphism $x: B\to X$
together with an associative (co)multiplication map $\mu:X\to X\otimes R$
satisfying the following properties:

(0) $\mu x=(\id_X\otimes 1) x$

(1) (pseudo-torsor) The multiplication $\mu$ together with $\id_X\otimes 1$ 
induce an isomorphism $X\otimes_BX\to X\otimes R$.

(2) (local triviality) The map $x$ is faithfully flat (this property does not
depend on the differential, see~\cite{app}).
\end{defn}

\subsubsection{}
\label{aff.def}
The definition above gives rise to a stack of groupoids on $\dga(k)$
in the topology generated by the faithfully flat maps. This is a (2-) functor
$\CC:\dga(k)\to \Grp$ such that for $B\in\dga(k)$ $\CC(B)$ is the groupoid
of $B$-torsors under $G$. 

For a fixed $B$-torsor $P$ one defines a fibred category
$\CC(P)$ over $\dgar(k)$ by the formula
$$\CC(P,A)=\{ A\otimes B\text{-torsors }\widetilde{P}
\text{ with a trivialization }\widetilde{P}\otimes_{A\otimes B}B\iso P\}$$

This is not yet the deformation functor we need since one cannot
expect that deformations with a dg base have no ``higher morphisms'' between 
them. Thus we define, for a given $B$-torsor $P$ and $A\in\dgar(k)$,
the following simplicial category $\CDD(P,A)$.
\begin{itemize}
\item{$\Ob\CDD(P,A)=\Ob\CC(P,A);$}
\item{$\uhom_{\CDD(P,A)}(x,y)_n=\Hom_{\CC(P_n,A)}(x_n,y_n).$}
\end{itemize}
Here $P_n$ is the torsor $\Omega_n\otimes P$ over $\Omega_n\otimes B$
and, similarly, $x_n, y_n$ are the torsors over $A\otimes\Omega_n\otimes B$.

Finally, we apply the simplicial nerve functor $\TN$ (see the Appendix)
to get a simplicial set from the simplicial category $\CDD(P,A)$:

$$ \ol{F}_P(A):=\TN(\CDD(P,A)).$$

\subsubsection{}
\label{non-aff}
Let now $S$ be an arbitrary scheme over $k$ and $P$ be a $S$-torsor under $G$.
We define then $F_P(A),\ A\in\dgar(k)$, by the formula
$$F_P(A)=\underset{\scriptsize{\begin{array}{l}
                       i:U\to S \\
                       U\text{ affine }
                  \end{array} }}{\holim}\ol{F}_{i^*P}(A),$$
the inverse homotopy limit being taken over all affine schemes over $S$.

As a result of our computations we will see in particular that 
$\widetilde{F}$ and $F$ are homotopy equivalent for affine $S$.

\subsubsection{}
\label{calc.triv}
Let us make some calculations. Suppose that $P$ is a trivial 
torsor over $B$.

Fix $A\in\dgar(k)$ and calculate the simplicial category $\CDD(A)$ --- we
omit $P$ from the notation since $P$ is supposed to be trivial.

The functor $\#$ forgetting the differential in dg objects, transforms
$B$-torsors to $B^{\#}$-torsors. 

\begin{lem}{triv-g}Let $P$ be a torsor over $B\otimes A$ trivial over $B$. 
Then $P^{\#}$ is a trivial $(B\otimes A)^{\#}$-torsor.
\end{lem}
\begin{pf}
Similarly to the usual algebraic geometry, one defines a cotangent complex
complex $L_{B/A}$ (see~\cite{ill}, ch.~II) for any graded (super) commutative 
$A$-algebra $B$ (since we work over a field of characteristic zero, we can 
use dg algebra resolutions of $B$ to calculate $L_{B/A}$). 
The cotangent complex commutes with the flat base change 
(see~\cite{ill}, II.2.2.3).

A map of graded (super) commutative algebras is {\em formally smooth} if it 
satisfies the left lifting property  with respect to surjective maps 
having a nilpotent kernel. As in~\cite{ill}, III, 3.1.2, a $A$-algebra $B$ 
whose cotangent complex is represented by a projective module, is formally 
smooth. 

Since the cotangent complex $L_{R/k}$ is finitely generated free, the base
change property and faithfully flat descent give that $L_{P/B}$ is finitely 
generated projective for any $B$-torsor $P$ under $G=\Spec(R)$. This implies 
that all graded torsors under $G$ are formally smooth. 

Then the map $P\to P\otimes_{B\otimes A}B\to B$ can be lifted to a map
$P\to B\otimes A$ splitting the structure map. This proves triviality of $P$.
\end{pf}

\subsubsection{}
\begin{cor}{triv-c}For a given $B\in\dga(k)$, the groupoid $\CC(A)$ is 
canonically (on $A$) equivalent to the following groupoid $\ol{\CC}(A)$
\begin{itemize}
   \item{ $\Ob\ol{\CC}(A)=\MC(\fm\otimes B\otimes\fg)$}
   \item{$\Hom_{\ol{\CC}(A)}(x,y)=\{g\in\exp((\fm\otimes B)^0\otimes\fg)| 
                                                             y=g(x)\}.$}
\end{itemize}
\end{cor}
\begin{pf}
According to~\ref{triv-g}, any $A\otimes B$-torsor trivial over $B$
has form $A\otimes B\otimes R$ as a graded algebra; its differential 
is defined by its restriction on $R$ which is a derivation
$\delta:R\to A\otimes B\otimes R$ trivialized by $A\to k$. This is given
by a Maurer-Cartan element of $\fm\otimes B\otimes\fg$. Any automorphism
of the graded torsor $(A\otimes B\otimes R)^{\#}$ is given by a 
$A\otimes B$-point of $G$. This one should map to the unit $B$-point
of $G$ under $A\to k$. This gives the second formula of the claim. 
\end{pf}

\subsubsection{}
\begin{prop}{triv-d}For a given commutative $k$-algebra $B$ the functor
$$F_{triv}:\dgar(k)\to\simpl$$
describing deformations of the trivial $B$-torsor under $G$, is naturally 
equivalent to the nerve of the Lie algebra $B\otimes\fg$.
\end{prop}
\begin{pf}
Let $B$ be a commutative $k$-algebra. According to~\ref{triv-c},
one has $\Ob\CDD(A)=\MC(\fm\otimes B\otimes\fg)$ is a singleton since 
$\fm\otimes B\otimes\fg$ is non-positively graded, so $\CDD(A)$
is actually a simplicial group. 

Re-denoting for simplicity $\CDD(A)$ by $\CDD$ and $\fm\otimes B\otimes\fg$ 
by $\fg$ we have
$$ \CDD_n=\{u\in\exp(\Omega_n\otimes\fg)^0|u(0)=0\}=\Stab_{G_n}(0)$$
in the notation of~\ref{somecalc}.

Now we have to find a natural equivalence from $\Sigma(\fg)$ to $\TN(\CDD)$.
Since $\MC(\fg)=\{0\}$, \Cor{pcs} furnishes us a principal fibration
with the base $\Sigma(\fg)$, the total space $G_{\bullet}$ with the group
$\Stab_{G_{\bullet}}(0)$ and a canonical pseudo-cross section.

The simplicial set $G_{\bullet}$ is isomorphic to 
$$ n\mapsto (\Omega_n\otimes\fg)^0$$
which is (a simplicial vector space and) a direct sum of simplicial
vector spaces of form $\Omega_{\bullet}^p$ which are all contractible by
\cite{le}, p.~44.

 Then Theorem 21.13 of~\cite{may} provides a canonical homotopy equivalence
$\Sigma(\fg)\to\ol{W}(\CDD)$ --- see Appendix. According to~\Lem{w=n},
this gives a canonical equivalence in question.
\end{pf}

\subsubsection{}Now, using the faithfully flat descent (see, e.g., ~\cite{app})
and taking into account the main result of~\cite{ddg}, we obtain the following

\begin{thm}{torsors}
Let $S$ be a scheme over a field $k$ of characteristic zero, $G$ be an affine
algebraic group and $P$ be a $S$-torsor under $G$. Let $\fg$ be the Lie
algebra of $G$ and $\fg_P$ be the induced by $P$ coherent sheaf of Lie 
algebras on $S$.

Then the formal stack of deformations of $P$ has form 
$$\CC(\Right\Gamma(S,\fg_P))$$
where $\Right\Gamma(S,\fg_P)$ is calculated using Thom-Sullivan $\Tot$ functor
of ~\cite{hdtc}. 
\end{thm}

\subsubsection{} A similar result holds in the context of~\cite{ka}
where $G$-local systems are considered.

For this one defines torsors under $G$ over a ``dg-ringed space'' $(X,\CO_X)$
where $X$ is a topological space endowed with a sheaf of commutative 
dg $k$-algebras $\CO_X$, using the topology generated by surjective open
covers of $X$. Then, given a good topological space $X$ and a torsor $P$
under $G$ over the ringed space $(X,k_X)$ ($\CO=k_X$ is the constant sheaf),
its deformations over a local dg artinian $k$-algebra $A$ are defined as 
$(X,A\otimes k_X)$-torsors under $G$ endowed with a trivialization over 
$(X,k_X)$. These deformations are governed locally by the sheaf of Lie 
algebras $\fg_P$ and, therefore, globally, by $\Right\Gamma(X,\fg_P)$.

This formula coincides with Kapranov's~\cite{ka}, 2.5.1. Our construction
of the formal moduli seems to be more ``honest'' then the {\em ad hoc }
definition 2.2.2 of~\cite{ka}: we start with a reasonable functor on 
artinian rings and look for a representing object. Of course, it is not 
absolutely ``honest'' since we knew the answer to obtain.

\section{Appendix: Nerve of a simplicial category}

Let $\CC$ be a simplicial category. We present here two versions of 
the nerve of $\CC$ --- $\TN(\CC)$ and $\ol{W}(\CC)$ --- and prove they are 
canonically homotopy equivalent under some assumptions which are fulfilled in 
the application we have in mind (see~\ref{moduli-G}).

As usual, for $x,y\in\Ob\CC$ we denote by $\uhom(x,y)$ the simplicial set
of arrows from $x$ to $y$.

\subsection{}One can consider $\CC$ as a simplicial object $\{\CC_n\}$ 
in the category $\Cat$ of categories; thus, applying the standard nerve 
construction to each one of the $\CC_n$, we get a bisimplicial set; 
its diagonal is denoted by $\TN(\CC)$.

More explicitly, $n$-simplices of $\TN(\CC)$ are the sequences
$$ [f_1|\ldots|f_n],\ f_i\in\uhom_n(v_{i-1},v_i)\text{ where }
v_0,\ldots, v_n\in\Ob\CC.$$
The faces and the degeneracies are given by the standard formulas.

\subsection{} Another construction is a minor generalization of the
$\ol{W}$ functor, see~\cite{may}, chapter IV. The $n$-simplices of 
$\ol{W}(\CC)$ are the sequences
$$ (g_1|\ldots|g_n),\ g_i\in\uhom_{n-i}(v_{i-1},v_i)\text{ where }
v_0,\ldots, v_n\in\Ob\CC.$$

The faces and the degeneracies are given by the following formulas
\begin{itemize}
\item{$d_i(g_1|\ldots|g_n)=(d_{i-1}g_1|\ldots|d_1g_{i-1}|g_{i+1}d_0g_i|
\ldots|g_n)$ for $i\not=0,n$ }
\item{$d_0(g_1|\ldots|g_n)=(g_2|\ldots|g_n)$} 
\item{$d_n(g_1|\ldots|g_n)=(d_{n-1}g_1|\ldots|d_1g_n)$}
\item{$s_0(g_1|\ldots|g_n)=(\id|g_1|\ldots|g_n)$} 
\item{$s_i(g_1|\ldots|g_n)=(s_{i-1}g_1|\ldots|s_0g_i|\id|g_{i+1}|\ldots|
                           g_n)$ for $i>0$.} 
\end{itemize}

\subsection{}One has the following canonical maps connecting $\TN(\CC)$ and
$\ol{W}(\CC)$
$$ \pi:\TN(\CC)\rlarrows\ol{W}(\CC):\rho$$
given by the formulas
$$ \pi[f_1|\ldots|f_n]=(d_0f_1|d_0^2f_2|\ldots|d_0^nf_n)$$
$$ \rho(g_1|\ldots|g_n)=[s_0g_1|s_0^2g_2|\ldots|s_0^ng_n].$$

\begin{lem}{w=n}Suppose that the simplicial category $\CC$ satisfies
the following properties

(a) for each $v,w\in\Ob\CC$ the simplicial set $\uhom(v,w)$ is Kan fibrant.

(b) for each $n$ the category $\CC_n$ is a groupoid.

Then the maps $\pi,\rho$ are homotopy equivalences.
\end{lem}
\begin{pf}
We will directly check that under the said assumptions the map $\pi$
satisfies the RLP with respect to the maps $\partial\Delta^n\to\Delta^n$.

A map $\Delta^n\to\ol{W}(\CC)$ is given by a collection
$$g=(g_1|\ldots|g_n),\ g_i\in\uhom_{n-i}(v_{i-1},v_i).$$
A map $\partial\Delta^n\to\TN(\CC)$ is given by a compatible collection
$$ x^i=[x^i_1|\ldots|x^i_{n-1}]\text{ with }\deg x^i_j=n-1,$$
the compatibility conditions being the conditions (A), (B) below.

(A) The condition $d_i(x^k)=d_{k-1}(x^i)$ for $i<k$ amounts to the system
\begin{itemize}
\item{$d_i(x_j^k)=d_{k-1}(x^i_j)\text{ for }j\leq i-1\text{ or }j\geq k+1$}
\item{$d_i(x_{i+1}^k\circ x_i^k)=d_{k-1}(x^i_i)\text{ for }i\leq k-2$}
\item{$d_i(x_{j+1}^k)=d_{k-1}(x^i_j)\text{ for }i<j\leq k-2$}
\item{$d_i(x_k^k)=d_{k-1}(x^i_k\circ x^i_{k-1})\text{ for }i\leq k-2$}
\item{$d_{k-1}(x^k_k\circ x^k_{k-1})=d_{k-1}(x^{k-1}_k\circ x^{k-1}_{k-1})$}
\end{itemize}

(B) The compatibilities of $x^i$ with $g$ say that $\pi(x^i)=d_i(g)$ which 
gives
\begin{equation*}
d_0^jx^i_j=\begin{cases}
              d_{i-j}g_j\text{ for }j\leq i-1 \\
              g_{i+1}\circ d_0g_i\text{ for }j= i\\
              g_{j+1}\text{ for } j>i.
           \end{cases} 
\end{equation*}

(C) Now, we have to construct a map $\Delta^n\to\TN(\CC)$ i.e. a collection
$[f_1|\ldots|f_n]$ satisfying the conditions $d_if=x^i,\ \pi(f)=g$ which
can be rewritten as a system
\begin{enumerate}
\item{$d_0^j(f_j)=g_j$}
\item{$d_i(f_j)=\begin{cases}
               x^i_j,\ \ i\geq j+1\\
               x^i_{j-1},\ i\leq j-2
              \end{cases}$}
\item{$ d_{j-1}(f_j)\circ d_{j-1}(f_{j-1})=x^{j-1}_{j-1}$}
\end{enumerate}

One is looking for $f_j$ by induction. For $j=1$ we have prescribed
values for $d_i(f_1)$, $i\not=1$. Thus, checking their compatibility
and using the condition (a) of the Lemma, we deduce that the system 
admits a solution $f_1$. 

Suppose that $f_i$ have already been found for $i<j$ so that the 
equations above are satisfied.
One checks first of all that for $j>1$ the equation (C.1) follows
from (C.2). Afterwards, one finds the value for $d_{j-1}(f_j)$ using
the condition (b) of the Lemma. Then we have the prescribed values
for $d_i(f_j)$ for all $i\not=j$ and we have only to check using the 
compatibility conditions (A), (B) that these prescriptions given by (C.2), 
(C.3), are compatible.
\end{pf}


\begin{thebibliography}{MMMM}

\bibitem[BG]{bg} Bousfield,~ Gugenheim, On PL de Rham theory and a rational
homotopy type, Memoirs of the A.M.S., vol.~8, {\bf \#179}(1976). 

\bibitem[D]{d} V.~Drinfeld, a letter to V.~Schechtman, September 1988.

\bibitem[GM1]{gm1} Goldman, Millson, The deformation theory of representations
   of fundamental groups of compact K\"ahler manifolds,
   Publ. Math. IHES, {\bf 67} (1988), 43--96. 

\bibitem[GM2]{gm2} Goldman, Millson, The homotopy invariance
of the Kuranishi space, Ill. J. Math., {\bf 34:2} (1990), 337--367

\bibitem[H1]{ddg} V.~Hinich, Descent of Deligne groupoids,
Intern. Math. Res. Notices, (1997), \# 5, 223--239.

\bibitem[H2]{haha} V.~Hinich, Homological algebra of homotopy algebras,
Comm. in algebra, {\bf 25(10)}(1997), 3291--3323.

\bibitem[H3]{app} V. Hinich, Rings with approximation property admit a
dualizing complex, Math. Nachrichten, {\bf 163}(1993), 289--296.

\bibitem[HS1]{hla} V.~Hinich, V.~Schechtman, Homotopy Lie algebras,
Adv. Soviet Math.,{\bf 16}(1993), part 2, 1--29.

\bibitem[HS2]{hdtc} V.~Hinich, V.~Schechtman, Deformation theory
and Lie algebra homology, Parts 1,2, Algebra Colloquium, {\bf 4}(1997),
pp. 213--240 and 291--316.

\bibitem[HS3]{hlha} V.~Hinich, V.~Schechtman, On the homotopy limit of 
homotopy algebras, Lecture Notes in Mathematics, {\bf 1289} (1987), 240--264.

\bibitem[HiSi]{hs} A.~Hirschowitz, C.~Simpson, Descente pour les $n$-champs,
Preprint {\tt math.AG/9807049}.

\bibitem[Ill]{ill} L.~Illusie, Complexe cotangent et d\'eformations I,
Lecture Notes in Math., {\bf 239} (1971).

\bibitem[Ka]{ka} M.~Kapranov, Injective resolutions of $BG$ and derived moduli 
               spaces of local systems, preprint {\tt alg-geom/9710027}

\bibitem[Ko]{ko} M.~Kontsevich, Topics in algebra. Deformation theory,
Lecture notes, 1994.

\bibitem[L]{le} D.~Lehmann, Th\'eorie homotopique des formes diff\'erentielles
(d'apr\`es D.~Sullivan), Ast\'erisque, {\bf 45}, 1977.

\bibitem[M]{may} J.~P.~May, Simplicial objects in algebraic topology, 1967.

\bibitem[Q1]{ha} D.~Quillen, Homotopy algebra, Lecture Notes in Math., 
{\bf 41} (1967).

\bibitem[Q2]{rht} D.~Quillen, Rational homotopy theory, Annals of Math., 
{\bf 90}(1969), 205--295.

\bibitem[Sc]{sc} M.~Schlessinger, Functors on Artin rings, 
Trans. A.M.S. {\bf 130}(1968), 208--222.

\end{thebibliography}
\end{document}